\documentclass[12pt]{amsart}
\usepackage{graphicx}
\usepackage{amsmath}
\usepackage{amssymb}
\usepackage{setspace}
\newtheorem{thm}{Theorem}[section]

\newtheorem{prop}[thm]{Proposition}
\theoremstyle{definition}
\newtheorem{defn}[thm]{Definition}
\theoremstyle{remark}
\newtheorem{rem}[thm]{Remark}
\theoremstyle{conclusion}

\numberwithin{equation}{section}

\begin{document}
\title[Multi-linear and multi-parameter pseudo-differential operators]{$L^{p}$ estimates for multi-linear and multi-parameter pseudo-differential operators }

\author{Wei Dai and Guozhen Lu}

\address{School of Mathematical Sciences, Beijing Normal University, Beijing 100875, P. R. China}
\email{daiwei@bnu.edu.cn}

\address{Department of Mathematics, Wayne State University, Detroit, MI 48202, U. S. A.}
\email{gzlu@math.wayne.edu}

\begin{abstract}
We establish the pseudo-differential variant of the $L^{p}$ estimates for multi-linear and multi-parameter Coifman-Meyer multiplier operators proved by C. Muscalu, J. Pipher, T. Tao and C. Thiele in \cite{MPTT1,MPTT2}.
\end{abstract}
\maketitle {\small {\bf Keywords:} Multi-linear and multi-parameter pseudo-differential operators; One-parameter and multi-parameter paraproducts; $L^{p}$ estimates; Coifman-Meyer theorem.\\

}

\section{Introduction}

\subsection{Background}
For $n\geq 1$ and $d\geq 1$, let $m$ be a bounded function in $\mathbb{R}^{nd}$, smooth away from the origin and satisfying H\"{o}rmander-Mikhlin conditions\footnote{Throughout this paper, $A\lesssim B$ means that there exists a universal constant $C>0$ such that $A\leq CB$. If necessary, we use explicitly $A\lesssim_{\star,\cdots,\star}B$ to indicate that there exists a positive constant $C_{\star,\cdots,\star}$ depending only on the quantities appearing in the subscript continuously such that $A\leq C_{\star,\cdots,\star}B$.\\ ${\bf *}$ Research of the second author is partly supported by a US NSF grant.}:
\begin{equation}\label{eq11}
    |\partial^{\alpha}m(\xi)|\lesssim\frac{1}{|\xi|^{|\alpha|}}
\end{equation}
for sufficiently many multi-indices $\alpha$. Denote by $T_{m}$ the $n$-linear operator defined by
\begin{equation}\label{eq12}
    T_{m}(f_{1},\cdots,f_{n})(x):=\int_{\mathbb{R}^{nd}}m(\xi)\hat{f_{1}}(\xi_{1})\cdots\hat{f_{n}}(\xi_{n})e^{2\pi ix\cdot(\xi_{1}+\cdots+\xi_{n})}d\xi,
\end{equation}
where $\xi=(\xi_{1},\cdots,\xi_{n})\in\mathbb{R}^{nd}$ and $f_{1},\cdots,f_{n}$ are Schwartz functions on $\mathbb{R}^{d}$. From the classical Coifman-Meyer theorem (see \cite{CJ, CM1,CM2,GT,KS}), we know that the operator $T_{m}$ extends to a bounded $n$-linear operator from $L^{p_{1}}(\mathbb{R}^{d})\times\cdots\times L^{p_{n}}(\mathbb{R}^{d})$ into $L^{p}(\mathbb{R}^{d})$, provided that $1<p_{1},\cdots,p_{n}\leq\infty$ and $\frac{1}{p}=\frac{1}{p_{1}}+\cdots+\frac{1}{p_{n}}>0$. When $n=2$, as a consequence of bilinear $T1$ theorem (see \cite{CJ, GT}), there is also a pseudo-differential variant of the classical Coifman-Meyer theorem for symbol $a\in BS^{0}_{1,0}(\mathbb{R}^{3d})$, that is, $a$ satisfies the differential inequalities
\begin{equation}\label{eq13}
    |\partial_{x}^{\gamma}\partial_{\xi}^{\alpha}\partial_{\eta}^{\beta}a(x,\xi,\eta)|\lesssim_{d,\alpha,\beta,\gamma}
    (1+|\xi|+|\eta|)^{-|\alpha|-|\beta|}
\end{equation}
for sufficiently many multi-indices $\alpha$, $\beta$, $\gamma$. Namely, let $T_{a}$ be the corresponding bilinear pseudo-differential operators defined by replacing $m$ with $a$ in \eqref{eq12}, then $T_{a}$ is bounded from $L^{p}(\mathbb{R}^{d})\times L^{q}(\mathbb{R}^{d})$ into $L^{r}(\mathbb{R}^{d})$, provided that $1<p, \, q\leq\infty$ and $\frac{1}{r}=\frac{1}{p}+\frac{1}{q}>0$ (see \cite{BT}, and see \cite{Be,Mu1,MS} for $d=1$ case). For large amounts of literature involving estimates for multi-linear Calder\'{o}n-Zygmund operators and multi-linear pseudo-differential operators, refer to e.g. \cite{BMNT2, CJ,CM2, GN, GT,GT1,KS,MS,MTT1}.

However, when we come into the situation that a differential operator (with different behaviors on different spatial variables $x_{i}$, $i=1,\cdots,d$) acts on a product of several functions (for instance, the bilinear form $\mathcal{D}_{1}^{\alpha}\mathcal{D}_{2}^{\beta}(fg)$, where $\widehat{\mathcal{D}_{1}^{\alpha}f}(\xi_{1},\xi_{2}):=|\xi_{1}|^{\alpha}\hat{f}(\xi_{1},\xi_{2})$ and $\widehat{\mathcal{D}_{2}^{\beta}f}(\xi_{1},\xi_{2}):=|\xi_{2}|^{\beta}\hat{f}(\xi_{1},\xi_{2})$ for $\alpha, \, \beta>0$),  we realize that  the necessity to investigate bilinear and bi-parameter operators $T_{m}^{(2)}$ defined by
\begin{equation}\label{eq14}
    T_{m}^{(2)}(f,g)(x):=\int_{\mathbb{R}^{4}}m(\xi,\eta)\hat{f}(\xi)\hat{g}(\eta)e^{2\pi ix\cdot(\xi+\eta)}d\xi d\eta,
\end{equation}
where the symbol $m$ is smooth away from the planes $(\xi_{1},\eta_{1})=(0,0)$ and $(\xi_{2},\eta_{2})=(0,0)$ in $\mathbb{R}^{2}\times\mathbb{R}^{2}$ and satisfying the less restrictive Marcinkiewicz conditions
\begin{equation}\label{eq15}
    |\partial_{\xi_{1}}^{\alpha_{1}}\partial_{\xi_{2}}^{\beta_{1}}\partial_{\eta_{1}}^{\alpha_{2}}\partial_{\eta_{2}}^{\beta_{2}}m(\xi,\eta)|
    \lesssim\frac{1}{|(\xi_{1},\eta_{1})|^{\alpha_{1}+\alpha_{2}}}\cdot\frac{1}{|(\xi_{2},\eta_{2})|^{\beta_{1}+\beta_{2}}}
\end{equation}
for sufficiently many multi-indices $\alpha=(\alpha_{1},\alpha_{2})$, $\beta=(\beta_{1},\beta_{2})$.  It becomes more complicated and difficult  to establish the $L^p$ estimates for  $T_{m}^{(2)}$ than in the one-parameter multilinear situations or  $L^p$ estimates for linear multi-parameter singular integrals (see e.g., \cite{FS} and \cite{J}). In \cite{MPTT1}, by using the duality lemma of $L^{p,\infty}$ presented in \cite{MTT1}, the $L^{1,\infty}$ sizes and energies technique developed in \cite{MTT2} and multi-linear interpolation (see e.g. \cite{GrTa,MTT2}), Muscalu, Pipher, Tao and Thiele proved the following $L^{p}$ estimates for $T_{m}^{(2)}$ (see also \cite{MS}, and for subsequent endpoint estimates see \cite{LM}).
\begin{thm}\label{22-multiplier}(\cite{MPTT1})
The bilinear operator $T_{m}^{(2)}$ defined by \eqref{eq14} maps $L^{p}(\mathbb{R}^{2})\times L^{q}(\mathbb{R}^{2})\rightarrow L^{r}(\mathbb{R}^{2})$ boundedly, as long as $1<p, \, q\leq\infty$ and $\frac{1}{r}=\frac{1}{p}+\frac{1}{q}>0$.
\end{thm}
In general, any collection of $n$ generic vectors $\xi_{1}=(\xi_{1}^{i})_{i=1}^{d}, \cdots, \xi_{n}=(\xi_{n}^{i})_{i=1}^{d}$ in $\mathbb{R}^{d}$ generates naturally the following collection of $d$ vectors in $\mathbb{R}^{n}$:
\begin{equation}\label{eq16}
    \bar{\xi}_{1}=(\xi^{1}_{j})_{j=1}^{n}, \,\,\,\,\,\, \bar{\xi}_{2}=(\xi^{2}_{j})_{j=1}^{n}, \,\,\,\,\, \cdots, \,\,\,\,\, \bar{\xi}_{d}=(\xi^{d}_{j})_{j=1}^{n}.
\end{equation}
Let $m=m(\xi)=m(\bar{\xi})$ be a bounded symbol in $L^{\infty}(\mathbb{R}^{dn})$ that is smooth away from the subspaces $\{\bar{\xi}_{1}=0\}\cup\cdots\cup\{\bar{\xi}_{d}=0\}$ and satisfying
\begin{equation}\label{eq17}
    |\partial^{\alpha_{1}}_{\bar{\xi}_{1}}\cdots\partial^{\alpha_{d}}_{\bar{\xi}_{d}}m(\bar{\xi})|\lesssim
    \prod_{i=1}^{d}|\bar{\xi_{i}}|^{-|\alpha_{i}|}
\end{equation}
for sufficiently many multi-indices $\alpha_{1},\cdots,\alpha_{d}$. Denote by $T_{m}^{(d)}$ the $n$-linear multiplier operator defined by
\begin{equation}\label{eq18}
    T_{m}^{(d)}(f_{1},\cdots,f_{n})(x):=\int_{\mathbb{R}^{dn}}m(\xi)\hat{f_{1}}(\xi_{1})\cdots\hat{f_{n}}(\xi_{n})
    e^{2\pi ix\cdot(\xi_{1}+\cdots+\xi_{n})}d\xi.
\end{equation}
In \cite{MPTT2}, Muscalu, Pipher, Tao and Thiele generalized Theorem \ref{22-multiplier} to the $n$-linear and $d$-parameter setting for any $n\geq1$, $d\geq2$, their result is stated in the following theorem (see also \cite{MS}).
\begin{thm}\label{dn-multiplier}(\cite{MPTT2})
For any $n\geq1$, $d\geq2$, the $n$-linear, $d$-parameter multiplier operator $T_{m}^{(d)}$ maps $L^{p_{1}}(\mathbb{R}^{d})\times\cdots\times L^{p_{n}}(\mathbb{R}^{d})\rightarrow L^{p}(\mathbb{R}^{d})$ boundedly, provided that $1<p_{1},\cdots,p_{n}\leq\infty$ and $\frac{1}{p}=\frac{1}{p_{1}}+\cdots+\frac{1}{p_{n}}>0$.
\end{thm}

\subsection{Main results}
The purpose of this paper is to prove the pseudo-differential variant of the $L^{p}$ estimates for multi-linear, multi-parameter Coifman-Meyer multiplier operators obtained in \cite{MPTT1,MPTT2} (see Theorem \ref{22-multiplier}, Theorem \ref{dn-multiplier}).

Suppose that $a(x,\xi,\eta)$ is a smooth symbol satisfying
\begin{equation}\label{eq19}
    |\partial_{x_{1}}^{\gamma_{1}}\partial_{x_{2}}^{\gamma_{2}}\partial_{\xi_{1}}^{\alpha_{1}}\partial_{\xi_{2}}^{\beta_{1}}
    \partial_{\eta_{1}}^{\alpha_{2}}\partial_{\eta_{2}}^{\beta_{2}}a(x,\xi,\eta)|\lesssim\frac{1}{(1+|(\xi_{1},\eta_{1})|)
    ^{\alpha_{1}+\alpha_{2}}}\cdot\frac{1}{(1+|(\xi_{2},\eta_{2})|)^{\beta_{1}+\beta_{2}}}
\end{equation}
for sufficiently many multi-indices $\alpha=(\alpha_{1},\alpha_{2})$, $\beta=(\beta_{1},\beta_{2})$, $\gamma=(\gamma_{1},\gamma_{2})$, and denote by $T_{a}^{(2)}$ the bilinear operator given by
\begin{equation}\label{eq110}
    T_{a}^{(2)}(f,g)(x):=\int_{\mathbb{R}^{4}}a(x,\xi,\eta)\hat{f}(\xi)\hat{g}(\eta)e^{2\pi ix\cdot(\xi+\eta)}d\xi d\eta.
\end{equation}
In this paper, we prove that the same $L^{p}$ estimates as $T_{m}^{(2)}$ in Theorem \ref{22-multiplier} hold true for the operator $T_{a}^{(2)}$. The main theorem of this article is the following result.
\begin{thm}\label{main}
The bilinear and bi-parameter pseudo-differential operator $T_{a}^{(2)}$ defined by \eqref{eq110} maps $L^{p}(\mathbb{R}^{2})\times L^{q}(\mathbb{R}^{2})\rightarrow L^{r}(\mathbb{R}^{2})$ boundedly, as long as $1<p, \, q\leq\infty$ and $\frac{1}{r}=\frac{1}{p}+\frac{1}{q}>0$.
\end{thm}

In general, let symbol $a=a(x,\xi)=a(x,\bar{\xi})\in C^{\infty}(\mathbb{R}^{d(n+1)})$ and satisfy the differential inequalities
\begin{equation}\label{eq111}
    |\partial_{x}^{\gamma}\partial_{\bar{\xi}_{1}}^{\alpha_{1}}\cdots\partial_{\bar{\xi}_{d}}^{\alpha_{d}}a(x,\bar{\xi})|
    \lesssim\prod_{i=1}^{d}(1+|\bar{\xi}_{i}|)^{-|\alpha_{i}|}
\end{equation}
for sufficiently many multi-indices $\alpha_{1},\cdots,\alpha_{d}$ and $\gamma$. Denote by $T_{a}^{(d)}$ the $n$-linear multiplier operator defined by
\begin{equation}\label{eq112}
    T_{a}^{(d)}(f_{1},\cdots,f_{n})(x):=\int_{\mathbb{R}^{dn}}a(x,\xi)\hat{f_{1}}(\xi_{1})\cdots\hat{f_{n}}(\xi_{n})
    e^{2\pi ix\cdot(\xi_{1}+\cdots+\xi_{n})}d\xi,
\end{equation}
then we can naturally generalize Theorem \ref{main} to the $n$-linear and $d$-parameter setting for any $n\geq1$, $d\geq2$ and obtain the pseudo-differential variant of Theorem \ref{dn-multiplier}. Our generalized theorem in this paper is the following.
\begin{thm}\label{generalized}
For any $n\geq1$, $d\geq2$, the $n$-linear, $d$-parameter pseudo-differential operator $T_{a}^{(d)}$ maps $L^{p_{1}}(\mathbb{R}^{d})\times\cdots\times L^{p_{n}}(\mathbb{R}^{d})\rightarrow L^{p}(\mathbb{R}^{d})$ boundedly, provided that $1<p_{1},\cdots,p_{n}\leq\infty$ and $\frac{1}{p}=\frac{1}{p_{1}}+\cdots+\frac{1}{p_{n}}>0$.
\end{thm}
\begin{rem}\label{generalization}
For simplicity, we will only prove Theorem \ref{main} (the bilinear and bi-parameter case, $n=d=2$) in this paper. However, it will be clear from the proof that we can extend the argument to the general $n$-linear, $d$-parameter setting (Theorem \ref{generalized}) straightforwardly.  
\end{rem}

\subsection{Outline of the proof strategy of our main theorems} In this subsection,  we would like to give an overview of our proof strategy of main theorems and indicate its additional difficulty and complexity compared with the case of $L^p$ estimates for one-parameter and multi-linear pseudo-differential operators of C. Muscalu \cite{Mu1, MS} and the $L^p$ estimates for multi-linear and  multi-parameter multiplier theorem of C. Muscalu, J. Pipher, T. Tao and C. Thiele \cite{MPTT1, MPTT2}.

By using the idea  by C. Muscalu in \cite{Mu1,MS} to prove the $L^{p}$ estimates for one-parameter (d=1) and bilinear pseudo-differential operators $T_{a}=T_{a}^{(1)}$, we will first show that the proof of Theorem \ref{main} can be essentially reduced to proving a localized variant of the bilinear and bi-parameter Coifman-Meyer theorem (Theorem \ref{22-multiplier}), that is, some kind of localized $L^{p}$ estimates of the localized bilinear and bi-parameter operator $T_{a}^{(2),(0,0,\vec{0})}$ given by
$$ T_{a}^{(2),(0,0,\vec{0})}(f,g)(x)=\int_{\mathbb{R}^{4}}m_{\vec{0}}(\xi,\eta)\hat{f}(\xi)\hat{g}(\eta)
  e^{2\pi ix\cdot(\xi+\eta)}d\xi d\eta\cdot\varphi_{0}^{'}\otimes\varphi_{0}^{''}(x).$$
Then, since the symbol $m_{\vec{0}}(\xi,\eta)$ of the operator $T_{a}^{(2),(0,0,\vec{0})}$ satisfies differential estimates \eqref{eq39} which is stronger than the Marcinkiewicz condition \eqref{eq15}, by making use of the inhomogeneous Littlewood-Paley dyadic decomposition \eqref{eq21} and Bony's paraproducts decomposition \eqref{eq22}, we can discretize the bilinear and bi-parameter operator $T_{a}^{(2),(0,0,\vec{0})}$ and reduce the proof of localized $L^{p}$ estimates of $T_{a}^{(2),(0,0,\vec{0})}$ to proving the localized $L^{p}$ estimates for discrete and localized bilinear and bi-parameter paraproduct operators of the form
$$\overrightarrow{\Pi}_{a,\mathcal{R}}^{(2),(0,0,\vec{0})}(f,g)(x)=\{\sum_{\substack{R=I\times J\in\mathcal{R}, \\ |I|,|J|\leq1}}c_{R}\frac{1}{|R|^{\frac{1}{2}}}\langle f,\varphi_{R}^{1}\rangle\langle g,\varphi_{R}^{2}\rangle\varphi_{R}^{3}(x)\}\cdot\varphi_{0}^{'}\otimes\varphi_{0}^{''}(x).$$
It's actually an inhomogeneous variant of the discretization procedure presented by Muscalu et al. in \cite{MPTT1} to prove the bi-parameter Coifman-Meyer theorem (Theorem \ref{22-multiplier}).

Now, in order to prove Theorem \ref{main}, we only have the task of proving localized $L^{p}$ estimates for the localized bilinear and bi-parameter paraproduct operator $\overrightarrow{\Pi}_{a,\mathcal{R}}^{(2),(0,0,\vec{0})}$ (Proposition \ref{paraproducts}). One can observe that the supports of functions $f, \, g$ and the dyadic rectangles $R=I\times J$ may get close to or far away from the integral region $R_{00}=I_{0}\times J_{0}$ in two different directions $x_{1}$ and $x_{2}$ due to the bi-parameter setting, thus the situations will be more complicated than the one-parameter case (d=1) considered in \cite{Be,Mu1,MS}. Since rapid decay factors can be derived from $\varphi_{R}^{3}$ when $R$ is sufficiently far away from $R_{00}$ in $x_{1}$ or $x_{2}$ directions (for example, $R\subseteq(5R_{00})^{c}$), we will split the localized bilinear and bi-parameter paraproduct operator $\overrightarrow{\Pi}_{a,\mathcal{R}}^{(2),(0,0,\vec{0})}$ into a summation of a ``main term", ``hybrid terms" and an ``error term" (see subsection 5.1).

Compared with the one-parameter case, there are mainly two key ingredients in our estimates of $\overrightarrow{\Pi}_{a,\mathcal{R}}^{(2),(0,0,\vec{0})}$ (see Section 5), one is the estimates of the ``main term" and ``hybrid terms", the other is the estimates of the discrete bilinear operators $\overrightarrow{\Pi}_{a,\mathcal{R}}^{(2),(0,0,\vec{0})}$ corresponding to bilinear operators involved in decomposition \eqref{paraproducts decompostion} which contain at least one of $\Pi^{1}_{ll}$ or $\Pi^{2}_{ll}$ in the tensor products, such as $\Pi^{1}_{lh}\otimes\Pi^{2}_{ll}$, $\Pi^{1}_{hl}\otimes\Pi^{2}_{ll}$, $\Pi^{1}_{hh}\otimes\Pi^{2}_{ll}$, $\Pi^{1}_{ll}\otimes\Pi^{2}_{ll}$, $\Pi^{1}_{ll}\otimes\Pi^{2}_{hh}$, $\Pi^{1}_{ll}\otimes\Pi^{2}_{hl}$ and $\Pi^{1}_{ll}\otimes\Pi^{2}_{lh}$, here we will only consider the case $\Pi^{1}_{ll}\otimes\Pi^{2}_{ll}$ without loss of generality. For the estimates of the ``main term" and ``hybrid terms" (see subsections 5.2 and 5.4), if the supports of $f, \, g$ are close to $R_{00}$ in both $x_{1}$ and $x_{2}$ directions ($supp \, f, \, supp \, g\subseteq15R_{00}$, say), we can apply the Coifman-Meyer theorem (Theorem \ref{22-multiplier}) or Theorem \ref{paraproduct estimates} directly; if for $i=1,2$, at least one of the supports of $f, \, g$ or dyadic rectangle $R$ is far away from $R_{00}$ in $x_{i}$ direction, we will obtain enough decay factors from $\langle f,\varphi_{R}^{1}\rangle\cdot\langle g,\varphi_{R}^{2}\rangle$ or $\varphi^{3}_{R}$; otherwise, if the supports of $f, \, g$ and dyadic rectangle $R$ are all close to $R_{00}$ in $x_{1}$ (or $x_{2}$) direction while at least one of the supports of $f, \, g$ are far away from $R_{00}$ in $x_{2}$ (or $x_{1}$) direction, we can apply the one-parameter paraproducts estimates (Theorem \ref{one-parameter paraproducts}) with respect to $x_{1}$ (or $x_{2}$) variable directly and obtain sufficient decay factors in $x_{2}$ (or $x_{1}$) direction so as to reach our conclusions. As to the estimates of the discrete bilinear operator $\overrightarrow{\Pi}_{a,\mathcal{R}}^{(2),(0,0,\vec{0})}$ corresponding to $\Pi^{1}_{ll}\otimes\Pi^{2}_{ll}$ (see subsection 5.5), one easily observe that at least two of the families of $L^{2}$-normalized bump functions $(\varphi_{I}^{i})_{I\in\mathcal{I}}$ for $i=1,2,3$ and two of $(\varphi_{J}^{j})_{J\in\mathcal{J}}$ for $j=1,2,3$ are nonlacunary respectively, which means that, when the supports of $f, \, g$ and dyadic rectangle $R$ are all close to $R_{00}$ in one direction (i.e. $I\subseteq5I_{0}$ or $J\subseteq5J_{0}$) but at least one of the supports of $f, \, g$ are far away from $R_{00}$ in the other direction, we won't be able to apply the one-parameter paraproducts estimates (Theorem \ref{one-parameter paraproducts}) with respect to $x_{1}$ or $x_{2}$ variable any more; however, we can take advantage of the additional properties that $|I|\sim|J|\sim1$ for every dyadic intervals $I\in\mathcal{I}$ and $J\in\mathcal{J}$ to obtain the convergence of both $\sum_{I\subseteq5I_{0}}$ and $\sum_{J\subseteq5J_{0}}$; the other parts of the estimates for the discrete bilinear operator $\overrightarrow{\Pi}_{a,\mathcal{R}}^{(2),(0,0,\vec{0})}$ corresponding to $\Pi^{1}_{ll}\otimes\Pi^{2}_{ll}$ are similar to the estimates of the standard discrete paraproduct operator corresponding to $\Pi^{1}_{lh}\otimes\Pi^{2}_{hl}$.

The rest of this paper is organized as follows. In Section 2 we give some useful notations and preliminary knowledge. In Section 3 we reduce the proof of Theorem \ref{main} to proving a localized variant of bilinear and bi-parameter Coifman-Meyer multiplier estimates (Proposition \ref{localized Coifman-Meyer}). Section 4 is devoted to reducing the proof of localized Coifman-Meyer multiplier estimates (Proposition \ref{localized Coifman-Meyer}) further to proving some localized discrete bilinear and bi-parameter paraproducts estimates (Proposition \ref{paraproducts}). In Section 5 we carry out the proof of Proposition \ref{paraproducts}, which completes the proof of our main theorem, Theorem \ref{main}.

\section{Notations and preliminaries}

Let $\varphi\in\mathcal{S}(\mathbb{R})$ be an even Schwartz function such that $supp \, \hat{\varphi}\subseteq[-\frac{4}{3},\frac{4}{3}]$ and $\hat{\varphi}(\xi)=1$ on $[-\frac{3}{4},\frac{3}{4}]$, and define $\psi\in\mathcal{S}(\mathbb{R})$ to be the Schwartz function whose Fourier transform satisfies $\hat{\psi}(\xi):=\hat{\varphi}(\frac{\xi}{2})-\hat{\varphi}(\xi)$ and $supp \, \hat{\psi}\subseteq[-\frac{8}{3},-\frac{3}{4}]\cup[\frac{3}{4},\frac{8}{3}]$, such that $0\leq\hat{\varphi}(\xi), \hat{\psi}(\xi)\leq1$. Then, for every integer $k\geq0$, we define $\widehat{\varphi_{k}}, \, \widehat{\psi_{k}}\in\mathcal{S}(\mathbb{R})$ by
\begin{equation*}
  \widehat{\varphi_{k}}(\xi):=\hat{\varphi}(\frac{\xi}{2^{k}}), \,\,\,\,\,\, \widehat{\psi_{k}}(\xi):=\hat{\psi}(\frac{\xi}{2^{k}})=\widehat{\varphi_{k+1}}(\xi)-\widehat{\varphi_{k}}(\xi)
\end{equation*}
and observe that
$$supp \, \widehat{\varphi_{k}}\subseteq[-\frac{4}{3}\cdot2^{k},\frac{4}{3}\cdot2^{k}], \,\,\,\,\,\, supp \, \widehat{\psi_{k}}\subseteq[-\frac{8}{3}\cdot2^{k},-\frac{3}{4}\cdot2^{k}]\cup[\frac{3}{4}\cdot2^{k},\frac{8}{3}\cdot2^{k}],$$
and $supp \, \widehat{\psi_{k}}\bigcap supp \, \widehat{\psi_{k'}}=\emptyset$ for any integers $k, \, k'\geq0$ such that $|k-k'|\geq2$, $supp \, \hat{\varphi}\bigcap supp \, \widehat{\psi_{k}}=\emptyset$ for any integer $k\geq1$.

We use the convention $\widehat{\psi_{-1}}(\xi):=\hat{\varphi}(\xi)$, then it is easy to see
\begin{equation}\label{eq20}
  1=\sum_{k\geq-1}\widehat{\psi_{k}}(\xi)
\end{equation}
for every $\xi\in\mathbb{R}$, as a consequence, one obtains the following inhomogeneous Littlewood-Paley dyadic decomposition of arbitrary functions $f, \, g\in\mathcal{S}'(\mathbb{R})$:
\begin{equation}\label{eq21}
  f=\sum_{k_{1}\geq-1}f\ast\psi_{k_{1}}, \,\,\,\,\,\, g=\sum_{k_{2}\geq-1}f\ast\psi_{k_{2}};
\end{equation}
furthermore, according to the criterion whether the support of a certain part of $\widehat{f\cdot g}$ contains the origin of momentum space $\mathbb{R}_{\xi}$ or not, we have Bony's paraproducts decomposition (see \cite{Bo,Tao,Taylor}) of the product $f\cdot g$:
\begin{eqnarray}\label{eq22}
 \nonumber  f\cdot g&=&\sum_{k_{1},k_{2}\geq-1}(f\ast\psi_{k_{1}})(g\ast\psi_{k_{2}}) \\
 \nonumber &=&\{\sum_{-1\leq k_{1}\leq k_{2}-2}+\sum_{-1\leq k_{2}\leq k_{1}-2}+\sum_{k_{1},k_{2}\geq-1, \, |k_{1}-k_{2}|\leq1}\}(f\ast\psi_{k_{1}})(g\ast\psi_{k_{2}})\\
  &=&\sum_{k\geq1}(f\ast\widetilde{\varphi_{k}})(g\ast\psi_{k})+\sum_{k\geq1}(f\ast\psi_{k})(g\ast\widetilde{\varphi_{k}})
  +\sum_{k\geq0}(f\ast\psi_{k})(g\ast\widetilde{\psi_{k}})\\
 \nonumber &&+\{(f\ast\varphi)(g\ast\psi)+(f\ast\psi)(g\ast\varphi)+(f\ast\varphi)(g\ast\varphi)\}\\
 \nonumber &=:&\Pi_{lh}(f,g)+\Pi_{hl}(f,g)+\Pi_{hh}(f,g)+\Pi_{ll}(f,g),
\end{eqnarray}
where $\widehat{\widetilde{\varphi_{k}}}(\xi):=\widehat{\varphi_{k-1}}(\xi)=\hat{\tilde{\varphi}}(\frac{\xi}{2^{k}})$ for any $k\geq1$, $\hat{\tilde{\varphi}}(\xi):=\hat{\varphi}(2\xi)$, and
$$\widetilde{\psi_{k}}:=\sum_{|k'-k|\leq1, \, k'\geq0}\psi_{k'}$$
for any $k\geq0$.

We use the notation $\langle\cdot,\cdot\rangle$ to denote the complex scalar $L^{2}$ inner product; and use $A^{c}$ to denote the complementary set of a set $A$.

An interval $I$ on the real line $\mathbb{R}$ is called dyadic if it is of the form $I=2^{-k}[n, \, n+1]$ for some $k, \, n\in\mathbb{Z}$, a rectangle $R$ on the plane $\mathbb{R}^{2}$ is called dyadic if there exist some dyadic intervals $I$, $J$ such that $R=I\times J$. Following \cite{MS}, we first give the following definitions.

\begin{defn}\label{bump functions}
For $J\subseteq\mathbb{R}$ an arbitrary interval, we say that a smooth function $\Phi_{J}$ is a bump adapted to $J$, if and only if the following inequalities hold:
\begin{equation}\label{eq23}
  |\Phi_{J}^{(l)}(x)|\lesssim_{l,\alpha}\frac{1}{|J|^{l}}\cdot\frac{1}{(1+\frac{dist(x,J)}{|J|})^{\alpha}}
\end{equation}
for every integer $\alpha\in\mathbb{N}$ and for sufficiently many derivatives $l\in\mathbb{N}$. If $\Phi_{J}$ is a bump adapted to $J$, we say that $|J|^{-\frac{1}{2}}\Phi_{J}$ is an $L^{2}$-normalized bump adapted to $J$.
\end{defn}

\begin{defn}\label{lacunary}
A family of $L^{2}$-normalized adapted bump functions $(\varphi_{I})_{I}$ is said to be nonlacunary if and only if for every $I$ one has
\[supp \, \widehat{\varphi_{I}}\subseteq[-4|I|^{-1}, \, 4|I|^{-1}].\]
A family of $L^{2}$-normalized adapted bump functions $(\varphi_{I})_{I}$ is said to be lacunary if and only if for any $I$ one has
\[supp \, \widehat{\varphi_{I}}\subseteq[-4|I|^{-1}, \, -\frac{1}{4}|I|^{-1}]\cup[\frac{1}{4}|I|^{-1}, \, 4|I|^{-1}].\]
\end{defn}

\begin{defn}
Let $\mathcal{I}$ be a finite set of dyadic intervals. A bilinear expression of the type
\begin{equation}\label{eq24}
  \Pi_{\mathcal{I}}(f,g)=\sum_{I\in\mathcal{I}}c_{I}\frac{1}{|I|^{\frac{1}{2}}}\langle f, \, \varphi_{I}^{1}\rangle\langle g, \, \varphi_{I}^{2}\rangle\varphi_{I}^{3}
\end{equation}
is called a bilinear discretized paraproduct if and only if $(c_{I})_{I}$ is a bounded sequence of complex numbers and at least two of the families of $L^{2}$-normalized bump functions $(\varphi_{I}^{j})_{I}$ for $j=1, \, 2, \, 3$ are lacunary in the sense of Definition \ref{lacunary}.
\end{defn}

Then the following is well-known, see e.g. \cite{LM,MS,MPTT1,MPTT2}.
\begin{thm}\label{one-parameter paraproducts}
Any bilinear discretized paraproduct $\Pi_{\mathcal{I}}$ has a bounded mapping $L^{p}\times L^{q}$ to $L^{r}$ as long as $1<p, \, q\leq\infty$ and $\frac{1}{r}=\frac{1}{p}+\frac{1}{q}>0$. Moreover, the implicit constants in the bounds depend only on $p,q,r$ and are independent of the cardinality of $\mathcal{I}$, provided that the sequence $(c_{I})_{I}$ in \eqref{eq24} is bounded by a universal constant.
\end{thm}

Consider two discretized classical paraproducts given by
\[\Pi_{1,\mathcal{I}}(f_{1}, \, g_{1})=\sum_{I\in\mathcal{I}}c_{I}\frac{1}{|I|^{\frac{1}{2}}}\langle f_{1}, \, \varphi_{I}^{1}\rangle\langle g_{1}, \, \varphi_{I}^{2}\rangle\varphi_{I}^{3}\]
and
\[\,\,\,\,\,\, \Pi_{2,\mathcal{J}}(f_{2}, \, g_{2})=\sum_{J\in\mathcal{J}}c_{J}\frac{1}{|J|^{\frac{1}{2}}}\langle f_{2}, \, \varphi_{J}^{1}\rangle\langle g_{2}, \, \varphi_{J}^{2}\rangle\varphi_{J}^{3},\]
and define the bi-parameter discretized paraproduct $\overrightarrow{\Pi}_{\mathcal{R}}$ by $\overrightarrow{\Pi}_{\mathcal{R}}=\Pi_{1,\mathcal{I}}\otimes\Pi_{2,\mathcal{J}}$ or, more generally, by
\begin{equation}\label{eq25}
  \overrightarrow{\Pi}_{\mathcal{R}}(f, \, g)=\sum_{R\in\mathcal{R}}c_{R}\frac{1}{|R|^{\frac{1}{2}}}\langle f, \, \varphi_{R}^{1}\rangle\langle g, \, \varphi_{R}^{2}\rangle\varphi_{R}^{3},
\end{equation}
where the numbers $c_{R}$ are all bounded, the sum is over dyadic rectangles of the form $R=I\times J$ and $\varphi_{R}^{j}$ is defined by $\varphi_{R}^{j}:=\varphi_{I}^{j}\otimes\varphi_{J}^{j}$ for $j=1, \, 2, \, 3$. We have the following $L^{p}$ estimates for bi-parameter discretized paraproduct $\overrightarrow{\Pi}_{\mathcal{R}}$ (for the proof, refer to \cite{LM,MS,MPTT1,MPTT2}).
\begin{thm}\label{paraproduct estimates}
Any discrete bi-parameter paraproduct \eqref{eq25} is bounded from $L^{p}(\mathbb{R}^{2})\times L^{q}(\mathbb{R}^{2})\rightarrow L^{r}(\mathbb{R}^{2})$ provided that $1<p, \, q\leq\infty$ and $\frac{1}{r}=\frac{1}{p}+\frac{1}{q}>0$.
\end{thm}
In general, the extension of Theorem \ref{paraproduct estimates} to the $n$-linear and $d$-parameter setting also holds (see \cite{LM,MPTT2}), that is, there is an analogue of the H\"{o}lder-type $L^{p}$ estimates stated in Theorem \ref{paraproduct estimates} for the discretized $d$-parameter paraproducts of the form
\begin{equation}\label{eq26}
  \overrightarrow{\Pi}_{\mathcal{R}}(f_{1},\cdots,f_{n})=\sum_{R\in\mathcal{R}}c_{R}\frac{1}{|R|^{\frac{n-1}{2}}}\langle f_{1}, \, \varphi_{R}^{1}\rangle\cdots\langle f_{n}, \, \varphi_{R}^{n}\rangle\varphi_{R}^{n+1},
\end{equation}
where the sum runs over the dyadic parallelepipeds $R=I_{1}\times\cdots\times I_{d}\subseteq\mathbb{R}^{d}$.

\section{Reduction to a localized variant of bilinear, bi-parameter Coifman-Meyer multiplier estimates}

In this section, by using the idea presented by C. Muscalu in \cite{Mu1,MS} to prove the $L^{p}$ estimates for one-parameter (d=1) and bilinear pseudo-differential operators $T_{a}=T_{a}^{(1)}$, we will show that the proof of our main result (Theorem \ref{main}), i.e. the H\"{o}lder-type $L^{p}$ estimates for operator $T_{a}^{(2)}$) can be essentially reduced to proving a localized variant of the Coifman-Meyer theorem (Theorem \ref{22-multiplier}).

To this end, we first pick two sequences of smooth functions $(\varphi_{n}^{'})_{n\in\mathbb{Z}}$ and $(\varphi_{m}^{''})_{m\in\mathbb{Z}}$ respectively, such that $supp \, \varphi_{n}^{'}\subseteq[n-1, \, n+1]$, $supp \, \varphi_{m}^{''}\subseteq[m-1, \, m+1]$ and
\begin{equation}\label{eq31}
  \sum_{n\in\mathbb{Z}}\varphi_{n}^{'}(x_{1})=\sum_{m\in\mathbb{Z}}\varphi_{m}^{''}(x_{2})=1
\end{equation}
for every $x=(x_{1}, \, x_{2})\in\mathbb{R}^{2}$. As a consequence, the bilinear and bi-parameter operator $T_{a}^{(2)}$ can be decomposed as follows
\begin{equation}\label{eq32}
  T_{a}^{(2)}=\sum_{n, \, m\in\mathbb{Z}}T_{a}^{(2), \, (n, \, m)},
\end{equation}
where $T_{a}^{(2), \, (n, \, m)}(f, \, g)(x):=T_{a}^{(2)}(f, \, g)(x)\cdot\varphi_{n}^{'}\otimes\varphi_{m}^{''}(x)$.

Now we claim that for every $n, \, m\in\mathbb{Z}$, one has estimates
\begin{equation}\label{eq33}
  \|T_{a}^{(2), \, (n, \, m)}(f,g)\|_{L^{r}(\mathbb{R}^{2})}\lesssim \|f\tilde{\chi}_{R_{nm}}\|_{L^{p}(\mathbb{R}^{2})}\cdot\|g\tilde{\chi}_{R_{nm}}\|_{L^{q}(\mathbb{R}^{2})}
\end{equation}
provided that $1<p, \, q\leq\infty$ and $\frac{1}{r}=\frac{1}{p}+\frac{1}{q}>0$, where the approximate cutoff functions $\tilde{\chi}_{R_{nm}}(x):=\tilde{\chi}_{I_{n}}\otimes\tilde{\chi}_{J_{m}}(x)$, $\tilde{\chi}_{I}(x_{1}):=(1+\frac{dist(x_{1},I)}{|I|})^{-100}$, $\tilde{\chi}_{J}(x_{2}):=(1+\frac{dist(x_{2},J)}{|J|})^{-100}$, rectangles $R_{nm}:=I_{n}\times J_{m}$, intervals $I_{n}:=[n-1,n+1]$, $J_{m}:=[m-1,m+1]$.

Suppose that we have proved this claim \eqref{eq33},  Theorem \ref{main} will also be proved as a corollary, because we have
\begin{eqnarray*}
  \|T_{a}^{(2)}(f,g)\|_{L^{r}}&\lesssim&(\sum_{n,m\in\mathbb{Z}}\|T_{a}^{(2),(n,m)}(f,g)\|_{L^{r}}^{r})^{\frac{1}{r}}
  \lesssim(\sum_{n,m\in\mathbb{Z}}\|f\tilde{\chi}_{R_{nm}}\|_{L^{p}}^{r}\|g\tilde{\chi}_{R_{nm}}\|_{L^{q}}^{r})^{\frac{1}{r}}\\
  &\lesssim&(\sum_{n,m\in\mathbb{Z}}\|f\tilde{\chi}_{R_{nm}}\|_{L^{p}}^{p})^{\frac{1}{p}}
  (\sum_{n,m\in\mathbb{Z}}\|g\tilde{\chi}_{R_{nm}}\|_{L^{q}}^{q})^{\frac{1}{q}}\lesssim\|f\|_{L^{p}}\cdot\|g\|_{L^{q}},
\end{eqnarray*}
as long as $1<p, \, q\leq\infty$ and $\frac{1}{r}=\frac{1}{p}+\frac{1}{q}>0$, where we have used the convergence of series $\sum_{k\geq1}k^{-s}$ for $s\gg1$ to obtain the last inequality. Therefore, we only have the task of proving the claim \eqref{eq33}.

Now arbitrarily fix $n_{0}, \, m_{0}\in\mathbb{Z}$, then the symbol of operator $T_{a}^{(2),(n_{0},m_{0})}$ can be written as
\begin{equation*}
  a(x,\xi,\eta)\varphi_{n_{0}}^{'}(x_{1})\varphi_{m_{0}}^{''}(x_{2})=a(x,\xi,\eta)\widetilde{\varphi_{n_{0}}^{'}}(x_{1})
  \varphi_{n_{0}}^{'}(x_{1})\widetilde{\varphi_{m_{0}}^{''}}(x_{2})\varphi_{m_{0}}^{''}(x_{2}),
\end{equation*}
where $\widetilde{\varphi_{n_{0}}^{'}}$, $\widetilde{\varphi_{m_{0}}^{''}}$ are smooth functions supported on the interval $\widetilde{I_{n_{0}}}:=[n_{0}-2,n_{0}+2]$, $\widetilde{J_{m_{0}}}:=[m_{0}-2,m_{0}+2]$ and that equal $1$ on the support of $\varphi_{n_{0}}^{'}$, $\varphi_{m_{0}}^{''}$, respectively. One can split the restricted symbol $a(x,\xi,\eta)\widetilde{\varphi_{n_{0}}^{'}}(x_{1})\widetilde{\varphi_{m_{0}}^{''}}(x_{2})$ as a Fourier series with respect to the $x$ variable and rewrite the symbol of $T_{a}^{(2),(n_{0},m_{0})}$ as
\begin{equation}\label{eq34}
  (\sum_{\vec{l}\in\mathbb{Z}^{2}}m_{\vec{l}}(\xi,\eta)e^{2\pi i(x_{1},x_{2})\cdot(l_{1},l_{2})})\cdot\varphi_{n_{0}}^{'}\otimes\varphi_{m_{0}}^{''}(x),
\end{equation}
where the Fourier coefficients
\begin{equation}\label{eq35}
  m_{\vec{l}}(\xi,\eta)=\int_{\mathbb{R}^{2}}a(x,\xi,\eta)\widetilde{\varphi_{n_{0}}^{'}}(x_{1})\widetilde{\varphi_{m_{0}}^{''}}(x_{2})e^{-2\pi ix\cdot\vec{l}}dx.
\end{equation}
Thus we have the decomposition
\begin{equation}\label{eq36}
  T_{a}^{(2),(n_{0},m_{0})}=\sum_{\vec{l}\in\mathbb{Z}^{2}}T_{a}^{(2),(n_{0},m_{0},\vec{l})},
\end{equation}
where
\begin{equation*}
  T_{a}^{(2),(n_{0},m_{0},\vec{l})}(f,g)(x)=\int_{\mathbb{R}^{4}}m_{\vec{l}}(\xi,\eta)e^{2\pi ix\cdot\vec{l}}\hat{f}(\xi)\hat{g}(\eta)
  e^{2\pi ix\cdot(\xi+\eta)}d\xi d\eta\cdot\varphi_{n_{0}}^{'}\otimes\varphi_{m_{0}}^{''}(x)
\end{equation*}
for every $\vec{l}\in\mathbb{Z}^{2}$. By applying invariant operator $L:=\frac{-\vec{l}\cdot\nabla_{x}}{2\pi i|\vec{l}|^{2}}$ with property
$L(e^{-2\pi ix\cdot\vec{l}})=e^{-2\pi ix\cdot\vec{l}}$ to the expression \eqref{eq35} of $m_{\vec{l}}(\xi,\eta)$ and integrating by parts sufficiently many times, we deduce from the estimates \eqref{eq19} of symbol $a(x,\xi,\eta)$ that
\begin{equation}\label{eq37}
  |\partial_{\xi_{1}}^{\alpha_{1}}\partial_{\xi_{2}}^{\beta_{1}}\partial_{\eta_{1}}^{\alpha_{2}}\partial_{\eta_{2}}^{\beta_{2}}m_{\vec{l}}(\xi,\eta)|
  \lesssim\frac{1}{(1+|\vec{l}|)^{M}}\cdot\frac{1}{(1+|(\xi_{1},\eta_{1})|)^{|\alpha|}}\cdot\frac{1}{(1+|(\xi_{2},\eta_{2})|)^{|\beta|}}
\end{equation}
for a sufficiently large number $M$ and sufficiently many multi-indices $\alpha=(\alpha_{1},\alpha_{2})$, $\beta=(\beta_{1},\beta_{2})$. One can observe from \eqref{eq37} that the Fourier coefficients $m_{\vec{l}}(\xi,\eta)$ decay rapidly in $|\vec{l}|$ away from the origin $\vec{0}$, which is acceptable for summation, it will be clear from our proof that we only need to consider the operator corresponding to $\vec{l}=\vec{0}$, which is given by
\begin{equation}\label{eq38}
  T_{a}^{(2),(n_{0},m_{0},\vec{0})}(f,g)(x)=\int_{\mathbb{R}^{4}}m_{\vec{0}}(\xi,\eta)\hat{f}(\xi)\hat{g}(\eta)
  e^{2\pi ix\cdot(\xi+\eta)}d\xi d\eta\cdot\varphi_{n_{0}}^{'}\otimes\varphi_{m_{0}}^{''}(x),
\end{equation}
where the symbol $m_{\vec{0}}(\xi,\eta)$ satisfies the following differential estimates
\begin{equation}\label{eq39}
  |\partial_{\xi_{1}}^{\alpha_{1}}\partial_{\xi_{2}}^{\beta_{1}}\partial_{\eta_{1}}^{\alpha_{2}}\partial_{\eta_{2}}^{\beta_{2}}m_{\vec{0}}(\xi,\eta)|
  \lesssim\frac{1}{(1+|(\xi_{1},\eta_{1})|)^{|\alpha|}}\cdot\frac{1}{(1+|(\xi_{2},\eta_{2})|)^{|\beta|}}
\end{equation}
for sufficiently many multi-indices $\alpha=(\alpha_{1},\alpha_{2})$, $\beta=(\beta_{1},\beta_{2})$.

Now assume that we have proved
\begin{equation}\label{eq310}
  \|T_{a}^{(2),(0,0,\vec{0})}(f,g)\|_{L^{r}(\mathbb{R}^{2})}\lesssim\|f\tilde{\chi}_{R_{00}}\|_{L^{p}(\mathbb{R}^{2})}
  \cdot\|g\tilde{\chi}_{R_{00}}\|_{L^{q}(\mathbb{R}^{2})},
\end{equation}
then we can infer from \eqref{eq38}, \eqref{eq310} and translation invariance that
\begin{eqnarray*}
  \|T_{a}^{(2),(n_{0},m_{0},\vec{0})}(f,g)\|_{L^{r}(\mathbb{R}^{2})}&=&
  \|T_{a}^{(2),(0,0,\vec{0})}(\tau_{n_{0}}^{x_{1}}\tau_{m_{0}}^{x_{2}}f,\tau_{n_{0}}^{x_{1}}\tau_{m_{0}}^{x_{2}}g)\|_{L^{r}(\mathbb{R}^{2})} \\
  &\lesssim&\|\tau_{n_{0}}^{x_{1}}\tau_{m_{0}}^{x_{2}}f\tilde{\chi}_{R_{00}}\|_{L^{p}(\mathbb{R}^{2})}
  \|\tau_{n_{0}}^{x_{1}}\tau_{m_{0}}^{x_{2}}g\tilde{\chi}_{R_{00}}\|_{L^{q}(\mathbb{R}^{2})} \\
  &=&\|f\tilde{\chi}_{R_{n_{0}m_{0}}}\|_{L^{p}(\mathbb{R}^{2})}\cdot\|g\tilde{\chi}_{R_{n_{0}m_{0}}}\|_{L^{q}(\mathbb{R}^{2})},
\end{eqnarray*}
where $\tau^{x_{1}}_{y}f(x):=f(x_{1}+y,x_{2})$ and $\tau^{x_{2}}_{y}f(x):=f(x_{1},x_{2}+y)$. Therefore, we can assume $n_{0}=m_{0}=0$, since $n_{0}, \, m_{0}\in\mathbb{Z}$ are chosen arbitrarily, we come to a conclusion that the proof of our claim \eqref{eq33}, or more precisely, the proof of Theorem \ref{main} can be reduced to proving a localized variant of the bilinear and bi-parameter Coifman-Meyer theorem, that is, the following proposition.
\begin{prop}\label{localized Coifman-Meyer}
Let bilinear operator $T_{a}^{(2),(0,0,\vec{0})}$ be defined by
\begin{equation}\label{eq311}
  T_{a}^{(2),(0,0,\vec{0})}(f,g)(x)=\int_{\mathbb{R}^{4}}m_{\vec{0}}(\xi,\eta)\hat{f}(\xi)\hat{g}(\eta)
  e^{2\pi ix\cdot(\xi+\eta)}d\xi d\eta\cdot\varphi_{0}^{'}\otimes\varphi_{0}^{''}(x),
\end{equation}
where the symbol $m_{\vec{0}}(\xi,\eta)$ satisfies differential estimates \eqref{eq39}, then we have
\begin{equation}\label{eq312}
   \|T_{a}^{(2),(0,0,\vec{0})}(f,g)\|_{L^{r}(\mathbb{R}^{2})}\lesssim\|f\tilde{\chi}_{R_{00}}\|_{L^{p}(\mathbb{R}^{2})}
  \cdot\|g\tilde{\chi}_{R_{00}}\|_{L^{q}(\mathbb{R}^{2})},
\end{equation}
provided that $1<p, \, q\leq\infty$ and $\frac{1}{r}=\frac{1}{p}+\frac{1}{q}>0$.
\end{prop}

\section{Discretization into localized bilinear and bi-parameter paraproducts estimates}

In this section, by making use of the inhomogeneous Littlewood-Paley decomposition \eqref{eq21} and Bony's paraproducts decomposition \eqref{eq22}, we will apply an inhomogeneous variant of the discretization procedure presented by Muscalu et al. in \cite{MPTT1} to reduce the operator $T_{a}^{(2),(0,0,\vec{0})}$ to averages of discrete and localized bilinear and bi-parameter paraproduct operators of the form
\begin{equation}\label{eq41}
  \overrightarrow{\Pi}_{a,\mathcal{R}}^{(2),(0,0,\vec{0})}(f,g)(x)=\{\sum_{\substack{R=I\times J\in\mathcal{R}, \\ |I|,|J|\leq1}}c_{R}\frac{1}{|R|^{\frac{1}{2}}}\langle f,\varphi_{R}^{1}\rangle\langle g,\varphi_{R}^{2}\rangle\varphi_{R}^{3}\}\cdot\varphi_{0}^{'}\otimes\varphi_{0}^{''}(x).
\end{equation}
Observe that the symbol $m_{\vec{0}}(\xi,\eta)$ of $T_{a}^{(2),(0,0,\vec{0})}$ doesn't have singularity near the planes $(\xi_{1},\eta_{1})=(0,0)$, $(\xi_{2},\eta_{2})=(0,0)$ in $\mathbb{R}^{2}\times\mathbb{R}^{2}$ and satisfies \eqref{eq39}, which is stronger than the Marcinkiewicz condition \eqref{eq15}, thus in the present case, it will be clear from the discretization procedure that one can use the $L^{1}$-normalized bump functions $\{\psi_{k}\}_{k\geq-1}$ (which are adapted to intervals of sizes $2^{-k}\lesssim1$ and of heights $2^{k}$) to carry out inhomogeneous Littlewood-Paley dyadic decomposition with respect to $x_{1}$ and $x_{2}$ variables respectively.  That is the reason why we can restrict further that the summation in \eqref{eq41} runs over dyadic intervals having the property that $|I|, \, |J|\lesssim1$.

We proceed the discretization procedure as follows. First, from the inhomogeneous one-parameter Littlewood-Paley decomposition \eqref{eq20}, we can see that the bilinear paraproducts decomposition \eqref{eq22} with respect to $x_{1}$ variable is equivalent to the following decomposition of symbol $1(\xi_{1},\eta_{1})$:
\begin{eqnarray}\label{eq42}
 \nonumber  1(\xi_{1},\eta_{1})&=& (\sum_{k_{1}\geq-1}\widehat{\psi_{k_{1}}}(\xi_{1}))(\sum_{k_{2}\geq-1}\widehat{\psi_{k_{2}}}(\eta_{1}))
 =\sum_{k_{1},k_{2}\geq-1}\widehat{\psi_{k_{1}}}(\xi_{1})\widehat{\psi_{k_{2}}}(\eta_{1}) \\
  &=&\sum_{k\geq1}\widehat{\widetilde{\varphi_{k}}}(\xi_{1})\widehat{\psi_{k}}(\eta_{1})+\sum_{k\geq1}\widehat{\psi_{k}}(\xi_{1})\widehat{\widetilde{\varphi_{k}}}(\eta_{1})
  +\sum_{k\geq0}\widehat{\psi_{k}}(\xi_{1})\widehat{\widetilde{\psi_{k}}}(\eta_{1})\\
 \nonumber &&+\{\hat{\varphi}(\xi_{1})\hat{\psi}(\eta_{1})+\hat{\psi}(\xi_{1})\hat{\varphi}(\eta_{1})+\hat{\varphi}(\xi_{1})\hat{\varphi}(\eta_{1})\}.
\end{eqnarray}
Similarly, we can decompose the symbol $1(\xi_{2},\eta_{2})$ as
\begin{eqnarray}\label{eq43}
  1(\xi_{2},\eta_{2})&=&\sum_{l\geq1}\widehat{\widetilde{\varphi_{l}}}(\xi_{2})\widehat{\psi_{l}}(\eta_{2})
  +\sum_{l\geq1}\widehat{\psi_{l}}(\xi_{2})\widehat{\widetilde{\varphi_{l}}}(\eta_{2})
  +\sum_{l\geq0}\widehat{\psi_{l}}(\xi_{2})\widehat{\widetilde{\psi_{l}}}(\eta_{2}) \\
 \nonumber &&+\{\hat{\varphi}(\xi_{2})\hat{\psi}(\eta_{2})+\hat{\psi}(\xi_{2})\hat{\varphi}(\eta_{2})+\hat{\varphi}(\xi_{2})\hat{\varphi}(\eta_{2})\}.
\end{eqnarray}
Note that $1(\xi,\eta)=1(\xi_{1},\eta_{1})\cdot1(\xi_{2},\eta_{2})$, one obtain immediately from \eqref{eq42}, \eqref{eq43} a decomposition of the symbol $1(\xi,\eta)$ as a sum of sixteen terms. We can also split the symbol $m_{\vec{0}}(\xi,\eta):=m_{\vec{0}}(\xi,\eta)\cdot1(\xi,\eta)$ as a sum of sixteen terms in the same way as $1(\xi,\eta)$, one of these terms is
\begin{equation*}
\sum_{k,l\geq1}m_{\vec{0}}(\xi,\eta)\widehat{\widetilde{\varphi_{k}}}(\xi_{1})\widehat{\psi_{k}}(\eta_{1})\widehat{\psi_{l}}(\xi_{2})\widehat{\widetilde{\varphi_{l}}}(\eta_{2})
  :=\sum_{k,l\geq1}m_{\vec{0}}(\xi,\eta)(\widehat{\widetilde{\varphi_{k}}}\otimes\widehat{\psi_{l}})(\xi)\cdot(\widehat{\psi_{k}}\otimes\widehat{\widetilde{\varphi_{l}}})(\eta).
\end{equation*}

Therefore, by splitting the symbol $m_{\vec{0}}(\xi,\eta)$ as above, one can decompose the operator $T_{a}^{(2),\vec{0}}$ given by
\begin{equation}\label{eq44}
  T_{a}^{(2),\vec{0}}(f,g)(x):=\int_{\mathbb{R}^{4}}m_{\vec{0}}(\xi,\eta)\hat{f}(\xi)\hat{g}(\eta)
  e^{2\pi ix\cdot(\xi+\eta)}d\xi d\eta
\end{equation}
into a sum of sixteen bilinear and bi-parameter paraproduct operators as follows:
\begin{spacing}{1.5}
\begin{equation}\label{paraproducts decompostion}
  \begin{array}{ll}
 \quad T_{a}^{(2), \, \vec{0}}(f, \, g)&\\
 =(\Pi^{1}_{lh}\otimes\Pi^{2}_{lh})(f,g)+(\Pi^{1}_{lh}\otimes\Pi^{2}_{hl})(f,g)+(\Pi^{1}_{lh}\otimes\Pi^{2}_{hh})(f,g)+(\Pi^{1}_{lh}\otimes\Pi^{2}_{ll})(f,g)&\\
 \quad+(\Pi^{1}_{hl}\otimes\Pi^{2}_{lh})(f,g)+(\Pi^{1}_{hl}\otimes\Pi^{2}_{hl})(f,g)+(\Pi^{1}_{hl}\otimes\Pi^{2}_{hh})(f,g)+(\Pi^{1}_{hl}\otimes\Pi^{2}_{ll})(f,g)&\\
 \quad+(\Pi^{1}_{hh}\otimes\Pi^{2}_{lh})(f,g)+(\Pi^{1}_{hh}\otimes\Pi^{2}_{hl})(f,g)+(\Pi^{1}_{hh}\otimes\Pi^{2}_{hh})(f,g)+(\Pi^{1}_{hh}\otimes\Pi^{2}_{ll})(f,g)&\\
 \quad+(\Pi^{1}_{ll}\otimes\Pi^{2}_{lh})(f,g)+(\Pi^{1}_{ll}\otimes\Pi^{2}_{hl})(f,g)+(\Pi^{1}_{ll}\otimes\Pi^{2}_{hh})(f,g)+(\Pi^{1}_{ll}\otimes\Pi^{2}_{ll})(f,g),&
\end{array}
\end{equation}
\end{spacing}
\setlength{\parindent}{0em}where $\Pi^{i}$ denotes one of the ``low-high", ``high-low", ``high-high" and ``low-low" paraproducts (defined in Section 2, \eqref{eq22}) with respect to $x_{i}$ variable for $i=1,2$, for instance, one of these operators can be expressed as
\begin{eqnarray}\label{eq45}
 &&T_{a,(lh,hl)}^{(2), \vec{0}}(f,g)(x):=(\Pi^{1}_{lh}\otimes\Pi^{2}_{hl})(f,g)(x)\\
 \nonumber &:=&\sum_{k,l\geq1}\int_{\mathbb{R}^{4}}m_{\vec{0}}(\xi,\eta)(f\ast(\widetilde{\varphi_{k}}\otimes\psi_{l}))^{\wedge}(\xi)
 (g\ast(\psi_{k}\otimes\widetilde{\varphi_{l}}))^{\wedge}(\eta)e^{2\pi ix\cdot(\xi+\eta)}d\xi d\eta.
\end{eqnarray}

\setlength{\parindent}{1.2em}
Observe that the nine operators in the decomposition \eqref{paraproducts decompostion} of $T_{a}^{(2),\vec{0}}$ (which don't contain the exponents $\Pi^{1}_{ll}$ or $\Pi^{2}_{ll}$ in the tensor products) are quite similar to the operator $T_{a,(lh,hl)}^{(2),\vec{0}}$, all of them can be reduced to averages of classical discrete bilinear paraproduct operators of the form \eqref{eq25} with restrictions $|I|,|J|\lesssim1$ and at least two of the families of $L^{2}$-normalized bump functions $(\varphi_{I}^{j})_{I}$ for $j=1,2,3$ are lacunary in the sense of Definition \ref{lacunary}, the same property also holds for $(\varphi_{J}^{j})_{J}$ ($j=1,2,3$).

But the situations are subtle for the other seven operators in the decomposition \eqref{paraproducts decompostion} of $T_{a}^{(2),\vec{0}}$ which contain at least one of components $\Pi^{1}_{ll}$ or $\Pi^{2}_{ll}$ in tensor products, such as $\Pi^{1}_{lh}\otimes\Pi^{2}_{ll}$, $\Pi^{1}_{hl}\otimes\Pi^{2}_{ll}$, $\Pi^{1}_{hh}\otimes\Pi^{2}_{ll}$, $\Pi^{1}_{ll}\otimes\Pi^{2}_{ll}$, $\Pi^{1}_{ll}\otimes\Pi^{2}_{hh}$, $\Pi^{1}_{ll}\otimes\Pi^{2}_{hl}$, $\Pi^{1}_{ll}\otimes\Pi^{2}_{lh}$. By the discretization procedure described below, one can reduce these seven operators to averages of discrete bilinear paraproduct operators of the form \eqref{eq25} with restrictions $|I|,|J|\lesssim1$ (at least one of $I$, $J$ satisfies $|I|\sim1$ or $|J|\sim1$), and for at least one of the two dyadic interval families $\mathcal{I}$ and $\mathcal{J}$ (here we assume the tensor product contains $\Pi^{1}_{ll}$ and hence suppose it is dyadic interval family $\mathcal{I}$ without loss of generality), one has $|I|\sim1$ for every $I\in\mathcal{I}$ and at least two of the families of $L^{2}$-normalized bump functions $(\varphi_{I}^{j})_{I\in\mathcal{I}}$ for $j=1,2,3$ are nonlacunary. Therefore, there are mainly two differences between these seven bilinear operators and the operator $T_{a,(lh,hl)}^{(2),\vec{0}}$. First, these seven operators can't be reduced to averages of classical discrete bilinear paraproduct operators of the form \eqref{eq25} which is applicable for Theorem \ref{paraproduct estimates}. Second, the one-parameter paraproducts estimates (Theorem \ref{one-parameter paraproducts}) can't be applied to each of the components $\Pi^{i}_{ll}$ ($i=1,2$) either. However, observe that $\Pi_{ll}^{i}$ are both summations of finite terms for $i=1,2$ and the dyadic intervals $I\in\mathcal{I}$ (corresponding to $\Pi^{1}_{ll}$), $J\in\mathcal{J}$ (corresponding to $\Pi^{2}_{ll}$) satisfy $|I|\sim1$ and $|J|\sim1$, so we can take advantage of the Coifman-Meyer theorem (Theorem \ref{22-multiplier}) and the fact that both $\sum_{I\subseteq5I_{0}}$ and $\sum_{J\subseteq5J_{0}}$ are finite summations to avoid the troubles of applying Theorem \ref{paraproduct estimates} and Theorem \ref{one-parameter paraproducts}(see subsection 5.5 in Section 5). It will be clear from the proof that the other parts of our arguments have nothing to do with the properties whether the families of $L^{2}$-normalized bump functions $(\varphi_{I}^{j})_{I\in\mathcal{I}}$ and $(\varphi_{J}^{j})_{J\in\mathcal{J}}$ for $j=1,2,3$ are lacunary or not (see subsection 5.2, 5.3 and 5.4 in Section 5), thus we can deal with these seven operators in a quite similar way as $T_{a,(lh,hl)}^{(2),\vec{0}}$.

In a word, we only need to consider the operator $T_{a,(lh,hl)}^{(2),\vec{0}}$ from now on, and the proof of Proposition \ref{localized Coifman-Meyer}, or more precisely, the proof of Theorem \ref{main} can be reduced to proving the following localized estimates for $T_{a,(lh,hl)}^{(2),\vec{0}}$:
\begin{equation}\label{eq47}
  \|T_{a,(lh,hl)}^{(2),\vec{0}}(f,g)\cdot\varphi_{0}^{'}\otimes\varphi_{0}^{''}\|_{L^{r}(\mathbb{R}^{2})}\lesssim
  \|f\widetilde{\chi}_{R_{00}}\|_{L^{p}(\mathbb{R}^{2})}\cdot\|g\widetilde{\chi}_{R_{00}}\|_{L^{q}(\mathbb{R}^{2})},
\end{equation}
as long as $1<p, \, q\leq\infty$ and $\frac{1}{r}=\frac{1}{p}+\frac{1}{q}>0$.

Now consider the trilinear form $\Lambda_{a,(lh,hl)}^{(2),\vec{0}}(f,g,h)$ associated to $T_{a,(lh,hl)}^{(2),\vec{0}}(f,g)$, which can be written as
\begin{eqnarray}\label{eq48}
   &&\Lambda_{a,(lh,hl)}^{(2),\vec{0}}(f,g,h):=\int_{\mathbb{R}^{2}}T_{a,(lh,hl)}^{(2),\vec{0}}(f,g)(x)h(x)dx\\
 \nonumber &=&\sum_{k,l\geq1}\int_{\xi+\eta+\gamma=0}m_{\vec{0},k,l}(\xi,\eta,\gamma)(f\ast(\widetilde{\varphi_{k}}\otimes\psi_{l}))^{\wedge}(\xi)
 (g\ast(\psi_{k}\otimes\widetilde{\varphi_{l}}))^{\wedge}(\eta)\quad\quad\quad\\
 \nonumber &&\quad\quad\quad\quad\quad\quad\quad\quad\quad\quad\quad\quad\quad\quad\quad\quad\quad\quad\quad\quad
 \times(h\ast(\psi'_{k}\otimes\psi'_{l}))^{\wedge}(\gamma)d\xi d\eta d\gamma,
\end{eqnarray}
where $\widehat{\psi'_{k}}(\gamma_{1}):=\widehat{\psi'}(\frac{\gamma_{1}}{2^{k}})$, $\widehat{\psi'_{l}}(\gamma_{2}):=\widehat{\psi'}(\frac{\gamma_{2}}{2^{l}})$ for any $k,l\in\mathbb{Z}$, $\psi^{'}$ is a Schwartz function such that $supp \, \widehat{\psi'}\subseteq[-4,-\frac{1}{16}]\cup[\frac{1}{16},4]$ and $\widehat{\psi'}=1$ on $[-\frac{10}{3},-\frac{1}{12}]\cup[\frac{1}{12},\frac{10}{3}]$, while $m_{\vec{0},k,l}(\xi,\eta,\gamma):=m_{\vec{0}}(\xi,\eta)\cdot(\lambda'_{k}\otimes\lambda^{''}_{l})(\xi,\eta,\gamma)$, where $\lambda'_{k}(\xi_{1},\eta_{1},\gamma_{1}):=\lambda'(\frac{\xi_{1}}{2^{k}},\frac{\eta_{1}}{2^{k}},\frac{\gamma_{1}}{2^{k}})$,
$\lambda^{''}_{l}(\xi_{2},\eta_{2},\gamma_{2}):=\lambda^{''}(\frac{\xi_{2}}{2^{l}},\frac{\eta_{2}}{2^{l}},\frac{\gamma_{2}}{2^{l}})$ for any $k,l\in\mathbb{Z}$, and $\lambda'\otimes\lambda^{''}$ is an appropriate smooth function supported on a slightly larger parallelepiped than $supp \, ((\widetilde{\varphi}\otimes\psi)^{\wedge}(\xi)(\psi\otimes\widetilde{\varphi})^{\wedge}(\eta)(\psi'\otimes\psi')^{\wedge}(\gamma))$, which equals $1$ on
$supp \, ((\widetilde{\varphi}\otimes\psi)^{\wedge}(\xi)(\psi\otimes\widetilde{\varphi})^{\wedge}(\eta)(\psi'\otimes\psi')^{\wedge}(\gamma))$. We can decompose $m_{\vec{0},k,l}(\xi,\eta,\gamma)$ as a Fourier series:
\begin{equation}\label{eq49}
  m_{\vec{0},k,l}(\xi,\eta,\gamma)=\sum_{\vec{n_{1}},\vec{n_{2}},\vec{n_{3}}\in\mathbb{Z}^{2}}C_{\vec{n_{1}},\vec{n_{2}},\vec{n_{3}}}^{k,l}
  e^{2\pi i(n^{'}_{1},n^{'}_{2},n^{'}_{3})\cdot(\xi_{1},\eta_{1},\gamma_{1})/2^{k}}
  e^{2\pi i(n^{''}_{1},n^{''}_{2},n^{''}_{3})\cdot(\xi_{2},\eta_{2},\gamma_{2})/2^{l}},
\end{equation}
where the Fourier coefficients $C_{\vec{n_{1}},\vec{n_{2}},\vec{n_{3}}}^{k,l}$($k,l\geq1$) are given by
\begin{equation}\label{eq410}
  C_{\vec{n_{1}},\vec{n_{2}},\vec{n_{3}}}^{k,l}=\int_{\mathbb{R}^{6}}m_{\vec{0},k,l}((2^{k}\xi_{1},2^{l}\xi_{2}),(2^{k}\eta_{1},2^{l}\eta_{2}),
  (2^{k}\gamma_{1},2^{l}\gamma_{2}))e^{-2\pi i(\vec{n_{1}}\cdot\xi+\vec{n_{2}}\cdot\eta+\vec{n_{3}}\cdot\gamma)}d\xi d\eta d\gamma.
\end{equation}
By taking advantage of the differential estimates \eqref{eq39} for symbol $m_{\vec{0}}(\xi,\eta)$, one deduce from \eqref{eq410} and integrating by parts sufficiently many times that
\begin{equation}\label{eq411}
  |C_{\vec{n_{1}},\vec{n_{2}},\vec{n_{3}}}^{k,l}|\lesssim \prod_{j=1}^{3}\frac{1}{(1+|\vec{n_{j}}|)^{M}}
\end{equation}
for any $k,l\geq1$, where $M$ is sufficiently large.

Then, by a straightforward calculation, we can rewrite \eqref{eq48} as
\begin{eqnarray}\label{eq412}
   &&\Lambda_{a,(lh,hl)}^{(2),\vec{0}}(f,g,h)=\sum_{\vec{n_{1}},\vec{n_{2}},\vec{n_{3}}\in\mathbb{Z}^{2},k,l\geq1}
   \Lambda_{\vec{n_{1}},\vec{n_{2}},\vec{n_{3}},k,l}(f,g,h)\\
 \nonumber &:=&\sum_{k,l\geq1}\sum_{\vec{n_{1}},\vec{n_{2}},\vec{n_{3}}\in\mathbb{Z}^{2}}C_{\vec{n_{1}},\vec{n_{2}},\vec{n_{3}}}^{k,l}\int_{\mathbb{R}^{2}}
  (f\ast(\widetilde{\varphi_{k}}\otimes\psi_{l}))(x-(2^{-k}n'_{1},2^{-l}n^{''}_{1}))\\
 \nonumber &&\times(g\ast(\psi_{k}\otimes\widetilde{\varphi_{l}}))(x-(2^{-k}n'_{2},2^{-l}n^{''}_{2}))
  (h\ast(\psi'_{k}\otimes\psi'_{l}))(x-(2^{-k}n'_{3},2^{-l}n^{''}_{3}))dx.
\end{eqnarray}
Since the rapid decay in \eqref{eq411} is acceptable for summation, we only need to consider the part of the trilinear form corresponding to $\vec{n_{1}}=\vec{n_{2}}=\vec{n_{3}}=(0,0)$:
\begin{equation}\label{eq413}
  \tilde{\Lambda}_{a,(lh,hl)}^{(2),\vec{0}}(f,g,h):=\sum_{k,l\geq1}\Lambda_{\vec{0},\vec{0},\vec{0},k,l}(f,g,h).
\end{equation}
By splitting the integral region $\mathbb{R}^{2}$ into the union of unit squares, the $L^{2}$-normalization procedure and simple calculations, we can rewrite \eqref{eq413} as
\begin{eqnarray*}
   &&\tilde{\Lambda}_{a,(lh,hl)}^{(2),\vec{0}}(f,g,h)\\
   &=&\sum_{k,l\geq1}\int_{0}^{1}\int_{0}^{1}\sum_{\substack{I \,\, dyadic, \\ |I|=2^{-k}}}\sum_{\substack{J \,\, dyadic, \\ |J|=2^{-l}}}\frac{1}{|I|^{\frac{1}{2}}}\cdot\frac{1}{|J|^{\frac{1}{2}}}\langle f,\varphi_{I}^{1,\nu'}\otimes\varphi_{J}^{1,\nu^{''}}\rangle
   \langle g,\varphi_{I}^{2,\nu'}\otimes\varphi_{J}^{2,\nu^{''}}\rangle \\
   &&\quad\quad\quad\quad\quad\quad\quad\quad\quad\quad\quad\quad\quad\quad\quad\quad\quad\quad\quad\quad\quad
   \times\langle h,\varphi_{I}^{3,\nu'}\otimes\varphi_{J}^{3,\nu^{''}}\rangle d\nu' d\nu^{''}\\
  &=:&\int_{0}^{1}\int_{0}^{1}\sum_{\substack{R=I\times J \,\, dyadic, \\ |I|,|J|<1}}\frac{1}{|R|^{\frac{1}{2}}}\langle f,\varphi_{R}^{1,\nu}\rangle\langle g,\varphi_{R}^{2,\nu}\rangle
  \langle h,\varphi_{R}^{3,\nu}\rangle d\nu,
\end{eqnarray*}
where $\varphi_{I}^{1,\nu'}$ is defined by $\varphi_{I}^{1,\nu'}(x_{1}):=2^{-\frac{k}{2}}\overline{\widetilde{\varphi_{k}}(2^{-k}(n+\nu')-x_{1})}$ and $I$ is the dyadic interval $[2^{-k}n,2^{-k}(n+1)]$, and all the $L^{2}$-normalized bump functions $\varphi_{I}^{j,\nu'}$ (adapted to $I$) and $\varphi_{J}^{j,\nu^{''}}$ (adapted to $J$) for $j=1,2,3$ can also be defined similarly in such a way respectively.

The bilinear operator corresponding to the trilinear form $\tilde{\Lambda}_{a,(lh,hl)}^{(2),\vec{0}}(f,g,h)$ can be written as
\begin{equation}\label{eq414}
  \widetilde{\overrightarrow{\Pi}}_{a}^{(2),\vec{0}}(f,g)(x)=\int_{0}^{1}\int_{0}^{1}\sum_{\substack{R=I\times J \,\, dyadic, \\ |I|,|J|<1}}\frac{1}{|R|^{\frac{1}{2}}}\langle f,\varphi_{R}^{1,\nu}\rangle\langle g,\varphi_{R}^{2,\nu}\rangle\varphi_{R}^{3,\nu}(x) d\nu.
\end{equation}
Since $\widetilde{\overrightarrow{\Pi}}_{a}^{(2),\vec{0}}(f,g)$ is an average of some discrete bilinear paraproduct operators depending on the parameters $\nu=(\nu_{1},\nu_{2})\in[0,1]^{2}$, it is enough to prove our localized estimates \eqref{eq312} in Proposition \ref{localized Coifman-Meyer} for each of them, uniformly with respect to parameters $\nu=(\nu_{1},\nu_{2})$. We will do this in the particular case when the parameters $\nu=(\nu_{1},\nu_{2})=(0,0)$, but the same argument works in general. In this case, we change our notation and rewrite the corresponding bilinear operator in \eqref{eq414} as
\begin{equation}\label{eq415}
 \overrightarrow{\Pi}_{a}^{(2),\vec{0}}(f,g)(x)=\sum_{\substack{R=I\times J \,\, dyadic, \\ |I|,|J|<1}}\frac{1}{|R|^{\frac{1}{2}}}\langle f,\varphi_{R}^{1}\rangle\langle g,\varphi_{R}^{2}\rangle\varphi_{R}^{3}(x).
\end{equation}

Now we reach a conclusion that in order to prove the localized Coifman-Meyer estimates in Proposition \ref{localized Coifman-Meyer}, we only need to prove that the bilinear operator $\overrightarrow{\Pi}_{a}^{(2),(0,0,\vec{0})}:=\overrightarrow{\Pi}_{a}^{(2),\vec{0}}\cdot\varphi_{0}^{'}\otimes\varphi_{0}^{''}$ satisfies estimates
\begin{equation}\label{eq416}
  \|\overrightarrow{\Pi}_{a}^{(2),(0,0,\vec{0})}(f,g)\|_{L^{r}(\mathbb{R}^{2})}\lesssim\|f\chi_{R_{00}}\|_{L^{p}(\mathbb{R}^{2})}
  \cdot\|g\chi_{R_{00}}\|_{L^{q}(\mathbb{R}^{2})}
\end{equation}
for any $1<p, \, q\leq\infty$ and $\frac{1}{r}=\frac{1}{p}+\frac{1}{q}>0$.

By Fatou's lemma, we can also restrict the summation in the definition of $\overrightarrow{\Pi}_{a}^{(2),(0,0,\vec{0})}$ on arbitrary finite set $\mathcal{R}$ of dyadic rectangles, and prove the estimates are unform with respect to different choices of the set $\mathcal{R}$. In a word, we have reduced the proof of Proposition \ref{localized Coifman-Meyer}, or more precisely, the proof of Theorem \ref{main} to proving the following localized estimates for discrete bilinear paraproducts $\overrightarrow{\Pi}_{a,\mathcal{R}}^{(2),(0,0,\vec{0})}$.
\begin{prop}\label{paraproducts}
Let localized and discrete bilinear paraproduct operator $\overrightarrow{\Pi}_{a,\mathcal{R}}^{(2),(0,0,\vec{0})}$ be defined by
\begin{equation}\label{eq417}
  \overrightarrow{\Pi}_{a,\mathcal{R}}^{(2),(0,0,\vec{0})}(f,g)(x)=\{\sum_{\substack{R=I\times J\in\mathcal{R}, \\ |I|,|J|\leq1}}c_{R}\frac{1}{|R|^{\frac{1}{2}}}\langle f,\varphi_{R}^{1}\rangle\langle g,\varphi_{R}^{2}\rangle\varphi_{R}^{3}\}\cdot\varphi_{0}^{'}\otimes\varphi_{0}^{''}(x),
\end{equation}
where $\mathcal{R}$ is an arbitrary finite set of dyadic rectangles and sequence $(c_{R})_{R}$ is bounded by a universal constant, then we have
\begin{equation}\label{eq418}
  \|\overrightarrow{\Pi}_{a,\mathcal{R}}^{(2),(0,0,\vec{0})}(f,g)\|_{L^{r}(\mathbb{R}^{2})}\lesssim\|f\chi_{R_{00}}\|_{L^{p}(\mathbb{R}^{2})}
  \cdot\|g\chi_{R_{00}}\|_{L^{q}(\mathbb{R}^{2})},
\end{equation}
as long as $1<p, \, q\leq\infty$ and $\frac{1}{r}=\frac{1}{p}+\frac{1}{q}>0$. Moreover, the implicit constants in the bounds depend only on $p, \, q, \,r$ and are independent of the cardinality of $\mathcal{R}$.
\end{prop}

\section{Proof of Theorem \ref{main}}

In this section, we prove our main result Theorem \ref{main} by carrying out the proof of Proposition \ref{paraproducts}.

\subsection{Strategy of the proof}
First observe the form of operator $\overrightarrow{\Pi}_{a,\mathcal{R}}^{(2),(0,0,\vec{0})}$. Since the integral region is $R_{00}:=I_{0}\times J_{0}:=[-1,1]\times[-1,1]$, we observe that for these terms that $R=I\times J\in\mathcal{R}$ is far away from $R_{00}$ in the summation in \eqref{eq417}, for instance, $R\subseteq(5R_{00})^{c}$, there will be rapid decay factors derived from the $L^{2}$-normalized bump function $\varphi_{R}^{3}$, and hence these terms will be small and easily estimated. Therefore, the main terms in the summation in \eqref{eq417} will be the ones that $R\subseteq5R_{00}$ (say), one easily deduce from the Coifman-Meyer theorem (Theorem \ref{22-multiplier}) or Theorem \ref{paraproduct estimates} that the main contribution in these cases comes from the cutoffs of $f$ and $g$ whose supports are not far away from $R_{00}$ (for instance, $f\chi_{15R_{00}}$ and $g\chi_{15R_{00}}$), since the function $\tilde{\chi}_{R_{00}}$ is bounded from below near the rectangle $R_{00}$; for other parts of $f$, $g$ which are supported far away from $R_{00}$, there will be rapid decay factors derived from $\langle f,\varphi_{R}^{1}\rangle\cdot\langle g,\varphi_{R}^{2}\rangle$, which is acceptable for summation.

According to the above analysis and observing that in the bi-parameter setting the dyadic rectangles $R=I\times J$ may get close to or far away from the integral region $R_{00}=I_{0}\times J_{0}$ in two different directions $x_{1}$ and $x_{2}$ which is more complicated than the one-parameter case, we split the bilinear paraproduct operator $\overrightarrow{\Pi}_{a,\mathcal{R}}^{(2),(0,0,\vec{0})}$ as follows:

 \begin{equation}\label{eq511}
  \overrightarrow{\Pi}_{a,\mathcal{R}}^{(2),(0,0,\vec{0})}:=\overrightarrow{\Pi}_{a,\mathcal{R},I}^{(2),(0,0,\vec{0})}
  +\overrightarrow{\Pi}_{a,\mathcal{R},II}^{(2),(0,0,\vec{0})}
  +\overrightarrow{\Pi}_{a,\mathcal{R},III}^{(2),(0,0,\vec{0})}+\overrightarrow{\Pi}_{a,\mathcal{R},IV}^{(2),(0,0,\vec{0})},
\end{equation}
where the main term
\begin{equation}\label{eq512}
  \overrightarrow{\Pi}_{a,\mathcal{R},I}^{(2),(0,0,\vec{0})}(f,g)(x):=\{\sum_{\substack{R\in\mathcal{R}, \, R\subseteq5R_{00}, \\ |I|,\,|J|\leq1}}\frac{c_{R}}{|R|^{\frac{1}{2}}}\langle f,\varphi_{R}^{1}\rangle\langle g,\varphi_{R}^{2}\rangle\varphi_{R}^{3}(x)\}\varphi_{0}^{'}\otimes\varphi_{0}^{''}(x),
\end{equation}
the error term
\begin{equation}\label{eq513}
   \overrightarrow{\Pi}_{a,\mathcal{R},II}^{(2),(0,0,\vec{0})}(f,g)(x):=\{\sum_{\substack{R\in\mathcal{R}, \, |I|,\,|J|\leq1, \\ I\subseteq(5I_{0})^{c},\,J\subseteq(5J_{0})^{c}}}\frac{c_{R}}{|R|^{\frac{1}{2}}}\langle f,\varphi_{R}^{1}\rangle\langle g,\varphi_{R}^{2}\rangle\varphi_{R}^{3}\}\varphi_{0}^{'}\otimes\varphi_{0}^{''}(x),
\end{equation}
and the hybrid terms
\begin{equation}\label{eq514}
   \overrightarrow{\Pi}_{a,\mathcal{R},III}^{(2),(0,0,\vec{0})}(f,g)(x):=\{\sum_{\substack{R\in\mathcal{R}, \, |I|,\,|J|\leq1, \\ I\subseteq5I_{0},\,J\subseteq(5J_{0})^{c}}}\frac{c_{R}}{|R|^{\frac{1}{2}}}\langle f,\varphi_{R}^{1}\rangle\langle g,\varphi_{R}^{2}\rangle\varphi_{R}^{3}(x)\}\varphi_{0}^{'}\otimes\varphi_{0}^{''}(x),
\end{equation}
\begin{equation}\label{eq515}
   \overrightarrow{\Pi}_{a,\mathcal{R},IV}^{(2),(0,0,\vec{0})}(f,g)(x):=\{\sum_{\substack{R\in\mathcal{R}, \, |I|,\,|J|\leq1, \\ I\subseteq(5I_{0})^{c},\,J\subseteq5J_{0}}}\frac{c_{R}}{|R|^{\frac{1}{2}}}\langle f,\varphi_{R}^{1}\rangle\langle g,\varphi_{R}^{2}\rangle\varphi_{R}^{3}(x)\}\varphi_{0}^{'}\otimes\varphi_{0}^{''}(x).
\end{equation}

\subsection{Estimates of the main term $\vec{\Pi}_{a,\mathcal{R},I}^{(2),(0,0,\vec{0})}$}
Observe that $R\subseteq5R_{00}$ in this situation, we can't obtain enough decay factors from $\varphi_{R}^{3}$ on the integral region $R_{00}$, but if one of $f$ or $g$ is supported far away from the region $R_{00}$, we will get decay factors from $\langle f,\varphi_{R}^{1}\rangle\cdot\langle g,\varphi_{R}^{2}\rangle$, which is acceptable for summation. To this end, let us decompose the functions $f$, $g$ as follows:
\begin{equation}\label{eq521}
  f=\sum_{n_{1},m_{1}\in\mathbb{Z}}f\chi_{\widetilde{I}_{n_{1}}}\chi_{\widetilde{J}_{m_{1}}}=:\sum_{n_{1},m_{1}\in\mathbb{Z}}f\chi_{\widetilde{R}_{n_{1}m_{1}}},
\end{equation}
\begin{equation}\label{eq522}
  g=\sum_{n_{2},m_{2}\in\mathbb{Z}}g\chi_{\widetilde{I}_{n_{2}}}\chi_{\widetilde{J}_{m_{2}}}=:\sum_{n_{2},m_{2}\in\mathbb{Z}}g\chi_{\widetilde{R}_{n_{2}m_{2}}},
\end{equation}
where $\widetilde{I}_{n_{i}}:=[n_{i},n_{i}+1)$, $\widetilde{J}_{m_{i}}:=[m_{i},m_{i}+1)$, $\widetilde{R}_{n_{i}m_{i}}:=\widetilde{I}_{n_{i}}\times\widetilde{J}_{m_{i}}$ for $i=1,2$. Now insert the two decompositions into the formula \eqref{eq512} for $\overrightarrow{\Pi}_{a,\mathcal{R},I}^{(2),(0,0,\vec{0})}$ and we get
\begin{equation}\label{eq523}
  \overrightarrow{\Pi}_{a,\mathcal{R},I}^{(2),(0,0,\vec{0})}(f,g)(x)=\sum_{n_{1},m_{1}\in\mathbb{Z}}\sum_{n_{2},m_{2}\in\mathbb{Z}}
  \overrightarrow{\Pi}_{a,\mathcal{R},I}^{(2),(0,0,\vec{0})}(f\chi_{\widetilde{R}_{n_{1}m_{1}}},g\chi_{\widetilde{R}_{n_{2}m_{2}}}).
\end{equation}

If all the $n_{1}$, $m_{1}$, $n_{2}$, $m_{2}$ are not far from zero, that is, $|n_{1}|, |m_{1}|, |n_{2}|, |m_{2}|\leq15$, then one gets from the Coifman-Meyer theorem (Theorem \ref{22-multiplier}), or more precisely, the discrete paraproducts estimates theorem (Theorem \ref{paraproduct estimates}) that
\begin{eqnarray}\label{eq524}
  &&\|\sum_{|n_{1}|,|m_{1}|\leq15}\sum_{|n_{2}|,|m_{2}|\leq15}\overrightarrow{\Pi}_{a,\mathcal{R},I}^{(2),(0,0,\vec{0})}
  (f\chi_{\widetilde{R}_{n_{1}m_{1}}},g\chi_{\widetilde{R}_{n_{2}m_{2}}})\|_{L^{r}}\\
\nonumber &\lesssim&\|f\sum_{|n_{1}|,|m_{1}|\leq15}\chi_{\widetilde{R}_{n_{1}m_{1}}}\|_{L^{p}}\cdot\|g\sum_{|n_{2}|,|m_{2}|\leq15}\chi_{\widetilde{R}_{n_{2}m_{2}}}\|_{L^{q}}
\lesssim\|f\tilde{\chi}_{R_{00}}\|_{L^{p}}\|g\tilde{\chi}_{R_{00}}\|_{L^{q}},
\end{eqnarray}
since $\tilde{\chi}_{R_{00}}$ is bounded from below near the rectangle $R_{00}$.

If, however, we faces the other fifteen different situations when at least one of $|n_{1}|$, $|m_{1}|$, $|n_{2}|$, $|m_{2}|$ is large, we mainly consider two kind of cases: first, there are at least one of $n_{1},n_{2}$ and one of $m_{1},m_{2}$ are large, for instance, suppose all of $|n_{1}|,|m_{1}|,|n_{2}|,|m_{2}|>15$ are large, then $\langle f\chi_{\widetilde{R}_{n_{1}m_{1}}},\varphi_{R}^{1}\rangle\cdot\langle g\chi_{\widetilde{R}_{n_{2}m_{2}}},\varphi_{R}^{2}\rangle$ will provide a decay factor of the type:
$$\frac{1}{(1+\frac{|n_{1}|-6}{|I|})^{N_{1}}}\cdot\frac{1}{(1+\frac{|m_{1}|-6}{|J|})^{M_{1}}}\cdot\frac{1}{(1+\frac{|n_{2}|-6}{|I|})^{N_{2}}}
\cdot\frac{1}{(1+\frac{|m_{2}|-6}{|J|})^{M_{2}}}$$
for sufficiently large number $N_{1}$, $M_{1}$, $N_{2}$ and $M_{2}$, which is acceptable for both the summations
$$\sum_{|n_{1}|,|n_{2}|>15}\sum_{I\subseteq5I_{0}} \,\,\,\,\,\,\,\,\,\,\,\, and \,\,\,\,\,\,\,\,\,\,\,\, \sum_{|m_{1}|,|m_{2}|>15}\sum_{J\subseteq5J_{0}}$$
on dyadic intervals $I$, $J$; second, there are at least one of $n_{1},n_{2}$ or at least one of $m_{1},m_{2}$ that is large, for instance, suppose $|n_{1}|>15$ is large, the other $m_{1}$, $n_{2}$, $m_{2}$ are not far from zero, then in a similar but simpler way as above, we deduce that $\langle f,\varphi_{R}^{1}\rangle$ will provide a decay factor of the type $(1+\frac{|n_{1}|-6}{|I|})^{-N_{1}}$ for $N_{1}$ sufficiently large, which is only enough for the summation $\sum_{|n_{1}|>15}\sum_{I\subseteq5I_{0}}$ on dyadic intervals $I$, we can apply the one-parameter paraproducts estimates (Theorem \ref{one-parameter paraproducts}) to solve the summation $\sum_{J\subseteq5J_{0}}$ on dyadic intervals $J$.

As analyzed above, we only consider the case that all of $|n_{1}|,|m_{1}|,|n_{2}|,|m_{2}|>15$ are large, the proofs of other cases are similar. For arbitrary $|n_{1}|,|m_{1}|,|n_{2}|,|m_{2}|>15$ and each fixed $R\subseteq5R_{00}$, since $\varphi_{R}^{j}=\varphi_{I}^{j}\otimes\varphi_{J}^{j}$ and $(\varphi_{I}^{j})_{I\in\mathcal{I}}$, $(\varphi_{J}^{j})_{J\in\mathcal{J}}$ are families of $L^{2}$-normalized bump functions adapted to intervals $I$, $J$ respectively for $j=1,2,3$, we deduce from H\"{o}lder's inequality the corresponding one-term bilinear operator satisfies the following estimates
\begin{eqnarray}\label{eq525}
&& \quad \|\frac{c_{R}}{|R|^{\frac{1}{2}}}\langle f\chi_{\widetilde{R}_{n_{1}m_{1}}},\varphi_{R}^{1}\rangle\langle g\chi_{\widetilde{R}_{n_{2}m_{2}}},\varphi_{R}^{2}\rangle\varphi_{R}^{3}\cdot\varphi_{0}^{'}\otimes\varphi_{0}^{''}\|_{L^{r}}\\
\nonumber &\lesssim&\frac{1}{|R|^{\frac{1}{2}}}(1+\frac{dist(\widetilde{I}_{n_{1}},I)}{|I|})^{-N_{1}}(1+\frac{dist(\widetilde{J}_{m_{1}},J)}{|J|})^{-M_{1}}
(\frac{1}{|R|^{\frac{1}{2}}}\|f\chi_{\widetilde{R}_{n_{1}m_{1}}}\|_{L^{p}}|R|^{\frac{p-1}{p}})
\cdot\\
\nonumber &&\times(1+\frac{dist(\widetilde{I}_{n_{2}},I)}{|I|})^{-N_{2}}(1+\frac{dist(\widetilde{J}_{m_{2}},J)}{|J|})^{-M_{2}}
(\frac{1}{|R|^{\frac{1}{2}}}\|g\chi_{\widetilde{R}_{n_{2}m_{2}}}\|_{L^{q}}|R|^{\frac{q-1}{q}})\frac{|R|^{\frac{1}{r}}}{|R|^{\frac{1}{2}}}\\
\nonumber &\lesssim&(1+\frac{dist(\widetilde{I}_{n_{1}},I)}{|I|})^{-N_{1}}(1+\frac{dist(\widetilde{J}_{m_{1}},J)}{|J|})^{-M_{1}}
(1+\frac{dist(\widetilde{I}_{n_{2}},I)}{|I|})^{-N_{2}}\\
\nonumber &&\quad\quad\quad\quad\quad\quad\quad\quad\quad\quad
\times(1+\frac{dist(\widetilde{J}_{m_{2}},J)}{|J|})^{-M_{2}}\|f\chi_{\widetilde{R}_{n_{1}m_{1}}}\|_{L^{p}}\cdot\|g\chi_{\widetilde{R}_{n_{2}m_{2}}}\|_{L^{q}}
\end{eqnarray}
for any $1<p, \, q\leq\infty$ and $\frac{1}{r}=\frac{1}{p}+\frac{1}{q}>0$, here we have used the facts that $$(1+\frac{dist(x_{1},I)}{|I|})^{N_{j}}|I|^{\frac{1}{2}}\varphi_{I}^{j}, \,\,\,\,\,\,\, (1+\frac{dist(x_{2},J)}{|J|})^{M_{j}}|J|^{\frac{1}{2}}\varphi_{J}^{j}$$
are also $L^{\infty}$-normalized bump functions adapted to dyadic intervals $I$, $J$ respectively for $j=1,2$, where $N_{1}$, $M_{1}$, $N_{2}$, $M_{2}$ are sufficiently large numbers (it will be enough for us to assume $N_{1},M_{1},N_{2},M_{2}\simeq1000$).

By using \eqref{eq512}, one can use the triangle inequality if $r\geq1$ and the subadditivity of $\|\cdot\|_{L^{r}}^{r}$ if $0<r<1$ to sum the contributions of every $R\subseteq5R_{00}$ with $|I|,|J|\leq1$ given by \eqref{eq525} together and obtain (we only present here the arguments for $0<r<1$, the cases $r\geq1$ can be treated similarly):
\begin{eqnarray}\label{eq526}
&&\quad\|\overrightarrow{\Pi}_{a,\mathcal{R},I}^{(2),(0,0,\vec{0})}(f\chi_{\widetilde{R}_{n_{1}m_{1}}},g\chi_{\widetilde{R}_{n_{2}m_{2}}})\|^{r}_{L^{r}}\\
\nonumber &\lesssim&\sum_{k,l\geq0}\sum_{\substack{R=I\times J\subseteq5R_{00}, \\ |I|=2^{-k}, |J|=2^{-l}}}[(1+\frac{|n_{1}|-6}{2^{-k}})^{-N_{1}}(1+\frac{|m_{1}|-6}{2^{-l}})^{-M_{1}}(1+\frac{|n_{2}|-6}{2^{-k}})^{-N_{2}}\\
\nonumber &&\quad\quad\quad\quad\quad\quad\quad\quad\quad\quad
\times(1+\frac{|m_{2}|-6}{2^{-l}})^{-M_{2}}]^{r}\|f\chi_{\widetilde{R}_{n_{1}m_{1}}}\|^{r}_{L^{p}}\cdot\|g\chi_{\widetilde{R}_{n_{2}m_{2}}}\|^{r}_{L^{q}}\\
 \nonumber &\lesssim&[\prod_{i=1,2}\frac{1}{(|n_{i}|-6)^{N_{i}}}\cdot\frac{1}{(|m_{i}|-6)^{M_{i}}}]^{r}
 \cdot\|f\chi_{\widetilde{R}_{n_{1}m_{1}}}\|^{r}_{L^{p}}\cdot\|g\chi_{\widetilde{R}_{n_{2}m_{2}}}\|^{r}_{L^{q}}
\end{eqnarray}
for any $|n_{1}|, |m_{1}|, |n_{2}|, |m_{2}|>15$.

One easily obtain that
\begin{equation}\label{eq527}
  (|n_{i}|-6)^{-\frac{N_{i}}{2}}\lesssim\min_{x_{1}\in\widetilde{I}_{n_{i}}}\tilde{\chi}_{I_{0}}(x_{1}), \quad\, (|m_{i}|-6)^{-\frac{M_{i}}{2}}\lesssim\min_{x_{2}\in\widetilde{J}_{m_{i}}}\tilde{\chi}_{J_{0}}(x_{2})
\end{equation}
for $i=1,2$ and every $|n_{1}|, |m_{1}|, |n_{2}|, |m_{2}|>15$, where $N_{1}, M_{1}, N_{2}, M_{2}\simeq1000$ are large enough.

Therefore, by using \eqref{eq526} and \eqref{eq527}, one can use the triangle inequality if $r\geq1$ and the subadditivity of $\|\cdot\|_{L^{r}}^{r}$ if $0<r<1$ to sum the contributions of $\overrightarrow{\Pi}_{a,\mathcal{R},I}^{(2),(0,0,\vec{0})}(f\chi_{\widetilde{R}_{n_{1}m_{1}}},g\chi_{\widetilde{R}_{n_{2}m_{2}}})$ together and obtain (we only present here the arguments for $r\geq1$, the cases $0<r<1$ can be treated similarly):
\begin{eqnarray}\label{eq528}
 &&\quad\,\|\sum_{|n_{1}|,|m_{1}|>15}\sum_{|n_{2}|,|m_{2}|>15}\overrightarrow{\Pi}_{a,\mathcal{R},I}^{(2),(0,0,\vec{0})}
 (f\chi_{\widetilde{R}_{n_{1}m_{1}}},g\chi_{\widetilde{R}_{n_{2}m_{2}}})\|_{L^{r}}\\
 \nonumber &\lesssim&\sum_{|n_{1}|,|m_{1}|>15}\sum_{|n_{2}|,|m_{2}|>15}\prod_{i=1,2}\frac{1}{(|n_{i}|-6)^{N_{i}}}\frac{1}{(|m_{i}|-6)^{M_{i}}}
 \|f\chi_{\widetilde{R}_{n_{1}m_{1}}}\|_{L^{p}}\|g\chi_{\widetilde{R}_{n_{2}m_{2}}}\|_{L^{q}}\\
 \nonumber &\lesssim&\sum_{|n_{1}|,|m_{1}|>15}\sum_{|n_{2}|,|m_{2}|>15}\prod_{i=1,2}\frac{1}{(|n_{i}|-6)^{\frac{N_{i}}{2}}}\frac{1}{(|m_{i}|-6)^{\frac{M_{i}}{2}}}
 \|f\tilde{\chi}_{R_{00}}\|_{L^{p}}\|g\tilde{\chi}_{R_{00}}\|_{L^{q}}\\
 \nonumber &\lesssim&\|f\tilde{\chi}_{R_{00}}\|_{L^{p}}\cdot\|g\tilde{\chi}_{R_{00}}\|_{L^{q}}.
\end{eqnarray}
Similar to estimates \eqref{eq524} and \eqref{eq528}, we can get the estimates for the other different fourteen cases, then we insert these estimates into the decomposition \eqref{eq523} and finally get the estimates of $\overrightarrow{\Pi}_{a,\mathcal{R},I}^{(2),(0,0,\vec{0})}$ as follows
\begin{eqnarray}\label{eq529}
 \nonumber  \|\overrightarrow{\Pi}_{a,\mathcal{R},I}^{(2),(0,0,\vec{0})}(f,g)\|_{L^{r}(\mathbb{R}^{2})}
   &=&\|\sum_{n_{1},m_{1}\in\mathbb{Z}}\sum_{n_{2},m_{2}\in\mathbb{Z}}
  \overrightarrow{\Pi}_{a,\mathcal{R},I}^{(2),(0,0,\vec{0})}(f\chi_{\widetilde{R}_{n_{1}m_{1}}},g\chi_{\widetilde{R}_{n_{2}m_{2}}})\|_{L^{r}}\\
  &\lesssim&\|f\tilde{\chi}_{R_{00}}\|_{L^{p}(\mathbb{R}^{2})}\cdot\|g\tilde{\chi}_{R_{00}}\|_{L^{q}(\mathbb{R}^{2})},
\end{eqnarray}
provided that $1<p, \, q\leq\infty$ and $\frac{1}{r}=\frac{1}{p}+\frac{1}{q}>0$, this concludes our estimates of the main term $\overrightarrow{\Pi}_{a,\mathcal{R},I}^{(2),(0,0,\vec{0})}$.

\subsection{Estimates of the error term $\vec{\Pi}_{a,\mathcal{R},II}^{(2),(0,0,\vec{0})}$}
Since $R=I\times J$ with $I\subseteq(5I_{0})^{c}$, $J\subseteq(5J_{0})^{c}$, $R$ is sufficiently far away from the integral region $R_{00}$, the operator $\overrightarrow{\Pi}_{a,\mathcal{R},II}^{(2),(0,0,\vec{0})}$ has sufficiently many rapid decay factors derived from $\varphi_{R}^{3}$ and can be considered as an error term.

One can decompose the operator $\overrightarrow{\Pi}_{a,\mathcal{R},II}^{(2),(0,0,\vec{0})}$ as
\begin{equation}\label{eq531}
  \overrightarrow{\Pi}_{a,\mathcal{R},II}^{(2),(0,0,\vec{0})}:=\sum_{|n|,|m|\geq5}\overrightarrow{\Pi}_{a,\mathcal{R}}^{nm},
\end{equation}
where
\begin{equation}\label{eq532}
  \overrightarrow{\Pi}_{a,\mathcal{R}}^{nm}(f,g)(x):=\{\sum_{\substack{R\in\mathcal{R}, \,R\subseteq R_{nm}, |I|,\,|J|\leq1, \\ I\subseteq(5I_{0})^{c},\,J\subseteq(5J_{0})^{c}}}\frac{c_{R}}{|R|^{\frac{1}{2}}}\langle f,\varphi_{R}^{1}\rangle\langle g,\varphi_{R}^{2}\rangle\varphi_{R}^{3}\}\varphi_{0}^{'}\otimes\varphi_{0}^{''}(x)
\end{equation}
for $|n|, \, |m|\geq5$. For arbitrary $|n|, \, |m|\geq5$ and each fixed $R\subseteq R_{nm}$, since $\varphi_{R}^{j}=\varphi_{I}^{j}\otimes\varphi_{J}^{j}$ and $(\varphi_{I}^{j})_{I\in\mathcal{I}}$, $(\varphi_{J}^{j})_{J\in\mathcal{J}}$ are families of $L^{2}$-normalized bump functions adapted to intervals $I$, $J$ respectively for $j=1,2,3$, we deduce from H\"{o}lder's inequality the corresponding one-term bilinear operator satisfies the following estimates
\begin{eqnarray}\label{eq533}
&&\|\frac{c_{R}}{|R|^{\frac{1}{2}}}\langle f,\varphi_{R}^{1}\rangle\langle g,\varphi_{R}^{2}\rangle\varphi_{R}^{3}\cdot\varphi_{0}^{'}\otimes\varphi_{0}^{''}\|_{L^{r}}\\
\nonumber &\lesssim&\frac{1}{|R|^{\frac{1}{2}}}(\frac{1}{|R|^{\frac{1}{2}}}\|f\tilde{\chi}_{R_{nm}}\|_{L^{p}}|R|^{\frac{p-1}{p}})
\cdot(\frac{1}{|R|^{\frac{1}{2}}}\|g\tilde{\chi}_{R_{nm}}\|_{L^{q}}|R|^{\frac{q-1}{q}})(1+\frac{dist(I,I_{0})}{|I|})^{-N}\\
\nonumber &&\quad\quad\quad\quad\quad\quad\quad\quad\quad\quad\quad\quad\quad\quad\quad\quad\quad\quad\quad\quad
\times(1+\frac{dist(J,J_{0})}{|J|})^{-M}|R|^{\frac{1}{r}-\frac{1}{2}}\\
\nonumber &\lesssim&(1+\frac{dist(I,I_{0})}{|I|})^{-N}(1+\frac{dist(J,J_{0})}{|J|})^{-M}\|f\tilde{\chi}_{R_{nm}}\|_{L^{p}}\cdot\|g\tilde{\chi}_{R_{nm}}\|_{L^{q}}
\end{eqnarray}
for any $1<p, \, q\leq\infty$ and $\frac{1}{r}=\frac{1}{p}+\frac{1}{q}>0$, here we have used the facts that $$(1+\frac{dist(x_{1},I)}{|I|})^{N}|I|^{\frac{1}{2}}\varphi_{I}^{j} \quad\quad\quad and \quad\quad\quad (1+\frac{dist(x_{2},J)}{|J|})^{M}|J|^{\frac{1}{2}}\varphi_{J}^{j}$$
are also $L^{\infty}$-normalized bump functions adapted to dyadic intervals $I$, $J$ respectively for $j=1,2,3$, where $N$, $M$ are sufficiently large numbers (it will be enough for us to assume $N, \, M\simeq1000$).

By using \eqref{eq532} and summing the contributions of every $R\subseteq R_{nm}$ given by \eqref{eq533}, we get the estimates of operator $\overrightarrow{\Pi}_{a,\mathcal{R}}^{nm}$ as follows:
\begin{eqnarray}\label{eq534}
&&\|\overrightarrow{\Pi}_{a,\mathcal{R}}^{nm}(f,g)\|_{L^{r}}\\
\nonumber &\lesssim&(\sum_{k,l\geq0}\sum_{\substack{R=I\times J\subseteq R_{nm}, \\ |I|=2^{-k}, |J|=2^{-l}}}(1+\frac{|n|-2}{2^{-k}})^{-N}(1+\frac{|m|-2}{2^{-l}})^{-M})\|f\tilde{\chi}_{R_{nm}}\|_{L^{p}}\|g\tilde{\chi}_{R_{nm}}\|_{L^{q}}\\
 \nonumber &\lesssim&\frac{1}{(|n|-2)^{N}}\cdot\frac{1}{(|m|-2)^{M}}\|f\tilde{\chi}_{R_{nm}}\|_{L^{p}}\cdot\|g\tilde{\chi}_{R_{nm}}\|_{L^{q}}
\end{eqnarray}
for any $|n|, \, |m|\geq5$ and $r\geq1$; if $0<r<1$, we can use the subadditivity of $\|\cdot\|_{L^{r}}^{r}$ to sum the contributions in a completely similar way and get the estimate \eqref{eq534}.

Since we have for arbitrary $|n|, \, |m|\geq5$,
\begin{eqnarray*}
&&(|n|-2)^{-200}(|m|-2)^{-200}\max_{x\in\mathbb{R}^{2}}\{(1+\frac{dist(x_{1},I_{n})}{|I_{n}|})^{-100}(1+\frac{dist(x_{2},J_{m})}{|J_{m}|})^{-100}\}\\
&\lesssim& \min_{x\in R_{nm}}\{(1+\frac{dist(x_{1},I_{0})}{|I_{0}|})^{-100}(1+\frac{dist(x_{2},J_{0})}{|J_{0}|})^{-100}\},
\end{eqnarray*}
and hence we infer that
\begin{equation}\label{eq535}
  (|n|-2)^{-\frac{N}{3}}(|m|-2)^{-\frac{M}{3}}|\tilde{\chi}_{R_{nm}}(x)|\lesssim|\tilde{\chi}_{R_{00}}(x)|
\end{equation}
for every $x\in\mathbb{R}^{2}$ and $|n|, \, |m|\geq5$, where $N, \, M\simeq1000$ are large enough.

Therefore, by using \eqref{eq531}, \eqref{eq534} and \eqref{eq535}, one can use the triangle inequality if $r\geq1$ and the subadditivity of $\|\cdot\|_{L^{r}}^{r}$ if $0<r<1$ to sum the contributions of $\overrightarrow{\Pi}_{a,\mathcal{R}}^{nm}$ together and obtain (we only present here the arguments for $0<r<1$, the cases $r\geq1$ can be treated similarly):
\begin{eqnarray*}
 &&\|\overrightarrow{\Pi}_{a,\mathcal{R},II}^{(2),(0,0,\vec{0})}(f,g)\|_{L^{r}}^{r}\lesssim\sum_{|n|,|m|\geq5}\|\overrightarrow{\Pi}_{a,\mathcal{R}}^{nm}\|_{L^{r}}^{r}\\
 \nonumber &\lesssim&\sum_{|n|,|m|\geq5}\frac{1}{(|n|-2)^{rN}}\cdot\frac{1}{(|m|-2)^{rM}}\|f\tilde{\chi}_{R_{nm}}\|_{L^{p}}^{r}\cdot\|g\tilde{\chi}_{R_{nm}}\|_{L^{q}}^{r}\\
 \nonumber &\lesssim&\sum_{|n|,|m|\geq5}(|n|-2)^{-\frac{N}{6}}(|m|-2)^{-\frac{M}{6}}\|f\tilde{\chi}_{R_{00}}\|_{L^{p}}^{r}
 \|g\tilde{\chi}_{R_{00}}\|_{L^{q}}^{r}\lesssim\|f\tilde{\chi}_{R_{00}}\|_{L^{p}}^{r}\|g\tilde{\chi}_{R_{00}}\|_{L^{q}}^{r},
\end{eqnarray*}
and hence we get the estimates for $\overrightarrow{\Pi}_{a,\mathcal{R},II}^{(2),(0,0,\vec{0})}$ as follows
\begin{equation}\label{eq536}
  \|\overrightarrow{\Pi}_{a,\mathcal{R},II}^{(2),(0,0,\vec{0})}(f,g)\|_{L^{r}(\mathbb{R}^{2})}\lesssim\|f\tilde{\chi}_{R_{00}}\|_{L^{p}(\mathbb{R}^{2})}
  \cdot\|g\tilde{\chi}_{R_{00}}\|_{L^{q}(\mathbb{R}^{2})},
\end{equation}
as long as $1<p, \, q\leq\infty$ and $\frac{1}{r}=\frac{1}{p}+\frac{1}{q}>0$, this concludes our estimates of the error term $\overrightarrow{\Pi}_{a,\mathcal{R},II}^{(2),(0,0,\vec{0})}$.

\subsection{Estimates of the hybrid terms $\vec{\Pi}_{a,\mathcal{R},III}^{(2),(0,0,\vec{0})}$ and $\vec{\Pi}_{a,\mathcal{R},IV}^{(2),(0,0,\vec{0})}$}
We will only estimate the hybrid term $\overrightarrow{\Pi}_{a,\mathcal{R},III}^{(2),(0,0,\vec{0})}$, since by symmetry the arguments for estimating $\overrightarrow{\Pi}_{a,\mathcal{R},IV}^{(2),(0,0,\vec{0})}$ is completely similar.

The operator $\overrightarrow{\Pi}_{a,\mathcal{R},III}^{(2),(0,0,\vec{0})}$ may be regarded as the ``hybrid" in two aspects. First, it behaves like the main term $\overrightarrow{\Pi}_{a,\mathcal{R},I}^{(2),(0,0,\vec{0})}$ in $x_{1}$ direction, because $I\subseteq5I_{0}$ and if one of the functions $f$, $g$ is supported far away from $I_{0}$ in $x_{1}$ direction, then $\langle f,\varphi_{I}^{1}\rangle$ or $\langle g,\varphi_{I}^{2}\rangle$ will provide sufficient decay factors. Second, similar to the error term $\overrightarrow{\Pi}_{a,\mathcal{R},II}^{(2),(0,0,\vec{0})}$, since $J\subseteq(5J_{0})^{c}$, decay factors can always be derived from $\varphi_{J}^{3}(x_{2})$, no matter whether the supports of $f$, $g$ are far away from $J_{0}$ in $x_{2}$ direction or not. Therefore, the proof strategy for the hybrid term $\overrightarrow{\Pi}_{a,\mathcal{R},III}^{(2),(0,0,\vec{0})}$ will be a reasonable combination of the arguments for main term and error term.

First, we split the functions $f$, $g$ only with respect to $x_{1}$ variable:
\begin{equation}\label{eq541}
  f=\sum_{n_{1}\in\mathbb{Z}}f\chi_{\widetilde{I}_{n_{1}}}, \,\,\,\,\,\,\,\,\,\,\, g=\sum_{n_{2}\in\mathbb{Z}}g\chi_{\widetilde{I}_{n_{2}}},
\end{equation}
and insert the two decompositions into the formula \eqref{eq514} for $\overrightarrow{\Pi}_{a,\mathcal{R},III}^{(2),(0,0,\vec{0})}$. We can decompose the operator $\overrightarrow{\Pi}_{a,\mathcal{R},III}^{(2),(0,0,\vec{0})}$ as
\begin{equation}\label{eq542}
  \overrightarrow{\Pi}_{a,\mathcal{R},III}^{(2),(0,0,\vec{0})}:=\sum_{n_{1},n_{2}\in\mathbb{Z}}\sum_{|m|\geq5}\overrightarrow{\Pi}_{a,\mathcal{R}}^{(n_{1},n_{2}),m},
\end{equation}
where
\begin{equation}\label{eq543}
  \overrightarrow{\Pi}_{a,\mathcal{R}}^{(n_{1},n_{2}),m}(f,g)(x):=\{\sum_{\substack{R\in\mathcal{R}, \, |I|, |J|\leq1, \\ I\subseteq5I_{0},\,J\subseteq(5J_{0})^{c}, \\ J\subseteq J_{m}}}\frac{c_{R}}{|R|^{\frac{1}{2}}}\langle f\chi_{\widetilde{I}_{n_{1}}},\varphi_{R}^{1}\rangle\langle g\chi_{\widetilde{I}_{n_{2}}},\varphi_{R}^{2}\rangle\varphi_{R}^{3}(x)\}\varphi_{0}^{'}\otimes\varphi_{0}^{''}(x)
\end{equation}
for every $n_{1},n_{2}\in\mathbb{Z}$ and $|m|\geq5$.

We mainly consider two kind of cases: first, at least one of $n_{1},n_{2}$ are large, there are three different subcases which can be estimated similarly, for instance, suppose both of $|n_{1}|,|n_{2}|>15$ are large, then $\langle f\chi_{\widetilde{I}_{n_{1}}},\varphi_{R}^{1}\rangle\cdot\langle g\chi_{\widetilde{I}_{n_{2}}},\varphi_{R}^{2}\rangle$ will provide a decay factor of the type:
$$\frac{1}{(1+\frac{|n_{1}|-6}{|I|})^{N_{1}}}\cdot\frac{1}{(1+\frac{|n_{2}|-6}{|I|})^{N_{2}}}$$
for sufficiently large number $N_{1}$ and $N_{2}$, which is acceptable for the summation $\sum_{|n_{1}|,|n_{2}|>15}\sum_{I\subseteq5I_{0}}$ on dyadic intervals $I$, at the same time, since $J\subseteq(5J_{0})^{c}$, $\varphi_{R}^{3}$ will provide sufficiently rapid decay factors on $|m|\geq5$ of the type
$$\frac{1}{(1+\frac{|m|-2}{|J|})^{M}}$$
for the summation $\sum_{|m|\geq5}\sum_{J\subseteq J_{m}}$; second, both of $n_{1}$, $n_{2}$ are not far from zero, we can apply directly the one-parameter paraproducts estimates (Theorem \ref{one-parameter paraproducts}) to solve the summation $\sum_{I\subseteq5I_{0}}$ on dyadic intervals $I$, then in a quite similar but simpler way as above, we deduce that $\varphi_{R}^{3}$ will provide enough decay factors on $|m|$ for the summation $\sum_{|m|\geq5}\sum_{J\subseteq J_{m}}$.

As analysed above, we only consider the case that both of $|n_{1}|,|n_{2}|>15$ are large, the proofs of other three cases are similar.  For arbitrary $|n_{1}|,|n_{2}|>15$, $|m|\geq5$ and each fixed $R=I\times J$ with $I\subseteq5I_{0}$, $J\subseteq J_{m}$, since $\varphi_{R}^{j}=\varphi_{I}^{j}\otimes\varphi_{J}^{j}$ and $(\varphi_{I}^{j})_{I\in\mathcal{I}}$, $(\varphi_{J}^{j})_{J\in\mathcal{J}}$ are families of $L^{2}$-normalized bump functions adapted to intervals $I$, $J$ respectively for $j=1,2,3$, we deduce from H\"{o}lder's inequality the corresponding one-term bilinear operator satisfies the following estimates
\begin{eqnarray}\label{eq544}
&& \,\,\, \|\frac{c_{R}}{|R|^{\frac{1}{2}}}\langle f\chi_{\widetilde{I}_{n_{1}}},\varphi_{R}^{1}\rangle\langle g\chi_{\widetilde{I}_{n_{2}}},\varphi_{R}^{2}\rangle\varphi_{R}^{3}\cdot\varphi_{0}^{'}\otimes\varphi_{0}^{''}\|_{L^{r}(\mathbb{R}^{2})}\\
\nonumber &\lesssim& (1+\frac{dist(\widetilde{I}_{n_{1}},I)}{|I|})^{-N_{1}}(1+\frac{dist(\widetilde{I}_{n_{2}},I)}{|I|})^{-N_{2}}\\
\nonumber && \quad\quad\quad\quad\quad\quad\quad\quad
\times\|\frac{1}{|J|^{\frac{1}{2}}}\langle\|f\chi_{\widetilde{I}_{n_{1}}}\|_{L_{x_{1}}^{p}},|\varphi_{J}^{1}|\rangle
\langle\|g\chi_{\widetilde{I}_{n_{2}}}\|_{L_{x_{1}}^{q}},|\varphi_{J}^{2}|\rangle\varphi_{J}^{3}\cdot\varphi_{0}^{''}\|_{L^{r}_{x_{2}}}\\
\nonumber &\lesssim& (1+\frac{dist(\widetilde{I}_{n_{1}},I)}{|I|})^{-N_{1}}(1+\frac{dist(\widetilde{I}_{n_{2}},I)}{|I|})^{-N_{2}}(1+\frac{dist(J,J_{0})}{|J|})^{-M}\\
\nonumber && \times\frac{1}{|J|^{\frac{1}{2}}}(\frac{1}{|J|^{\frac{1}{2}}}\|f\chi_{\widetilde{I}_{n_{1}}}\tilde{\chi}_{J_{m}}\|_{L^{p}(\mathbb{R}^{2})}|J|^{1-\frac{1}{p}})
(\frac{1}{|J|^{\frac{1}{2}}}\|g\chi_{\widetilde{I}_{n_{2}}}\tilde{\chi}_{J_{m}}\|_{L^{q}(\mathbb{R}^{2})}|J|^{1-\frac{1}{q}})|J|^{\frac{1}{r}-\frac{1}{2}}\\
\nonumber &\lesssim& (1+\frac{dist(\widetilde{I}_{n_{1}},I)}{|I|})^{-N_{1}}(1+\frac{dist(\widetilde{I}_{n_{2}},I)}{|I|})^{-N_{2}}(1+\frac{dist(J,J_{0})}{|J|})^{-M}\\
\nonumber && \quad\quad\quad\quad\quad\quad\quad\quad\quad\quad\quad\quad\quad\quad\quad
\times\|f\chi_{\widetilde{I}_{n_{1}}}\tilde{\chi}_{J_{m}}\|_{L^{p}(\mathbb{R}^{2})}\cdot\|g\chi_{\widetilde{I}_{n_{2}}}\tilde{\chi}_{J_{m}}\|_{L^{q}(\mathbb{R}^{2})}
\end{eqnarray}
for any $1<p, \, q\leq\infty$ and $\frac{1}{r}=\frac{1}{p}+\frac{1}{q}>0$, here we have used the facts that $(1+\frac{dist(x_{1},I)}{|I|})^{N_{j}}|I|^{\frac{1}{2}}\varphi_{I}^{j}$ is also an $L^{\infty}$-normalized bump function adapted to dyadic interval $I$ for $j=1,2$ and $(1+\frac{dist(x_{2},J)}{|J|})^{M}|J|^{\frac{1}{2}}\varphi_{J}^{3}$ is also an $L^{\infty}$-normalized bump function adapted to dyadic interval $J$, where $N_{1}$, $N_{2}$, $M$ are sufficiently large numbers (it will be enough for us to assume $N_{1},N_{2},M\simeq1000$).

By using \eqref{eq543}, one can use the triangle inequality if $r\geq1$ and the subadditivity of $\|\cdot\|_{L^{r}}^{r}$ if $0<r<1$ to sum the contributions of every $R=I\times J$ with $I\subseteq5I_{0}$, $J\subseteq J_{m}$ and $|I|,|J|\leq1$ given by \eqref{eq544} together and obtain (we only present here the arguments for $0<r<1$, the cases $r\geq1$ can be treated similarly):
\begin{eqnarray}\label{eq545}
&&\quad\|\overrightarrow{\Pi}_{a,\mathcal{R}}^{(n_{1},n_{2}),m}(f,g)\|^{r}_{L^{r}(\mathbb{R}^{2})}\\
\nonumber &\lesssim&\sum_{k,l\geq0}\sum_{\substack{I\subseteq5I_{0}, \\ |I|=2^{-k}}}\sum_{\substack{J\subseteq J_{m}, \\ |J|=2^{-l}}}[(1+\frac{|n_{1}|-6}{2^{-k}})^{-N_{1}}(1+\frac{|n_{2}|-6}{2^{-k}})^{-N_{2}}(1+\frac{|m|-2}{2^{-l}})^{-M}]^{r}\\
\nonumber &&\quad\quad\quad\quad\quad\quad\quad\quad\quad\quad\quad\quad\quad\quad\quad
\times\|f\chi_{\widetilde{I}_{n_{1}}}\tilde{\chi}_{J_{m}}\|^{r}_{L^{p}(\mathbb{R}^{2})}
\cdot\|g\chi_{\widetilde{I}_{n_{2}}}\tilde{\chi}_{J_{m}}\|^{r}_{L^{q}(\mathbb{R}^{2})}\\
 \nonumber &\lesssim&[\prod_{i=1,2}\frac{1}{(|n_{i}|-6)^{N_{i}}}]^{r}\frac{1}{(|m|-2)^{Mr}}
 \cdot\|f\chi_{\widetilde{I}_{n_{1}}}\tilde{\chi}_{J_{m}}\|^{r}_{L^{p}(\mathbb{R}^{2})}
 \cdot\|g\chi_{\widetilde{I}_{n_{2}}}\tilde{\chi}_{J_{m}}\|^{r}_{L^{q}(\mathbb{R}^{2})}
\end{eqnarray}
for any $|n_{1}|, |n_{2}|>15$ and $|m|\geq5$.

Since we have for arbitrary $|m|\geq5$,
\begin{equation*}
(|m|-2)^{-200}\max_{x_{2}\in\mathbb{R}}(1+\frac{dist(x_{2},J_{m})}{|J_{m}|})^{-100}\lesssim\min_{x_{2}\in J_{m}}
(1+\frac{dist(x_{2},J_{0})}{|J_{0}|})^{-100},
\end{equation*}
and hence we infer that
\begin{equation}\label{eq546}
  (|m|-2)^{-\frac{M}{4}}|\tilde{\chi}_{J_{m}}(x_{2})|\lesssim|\tilde{\chi}_{J_{0}}(x_{2})|
\end{equation}
for every $x_{2}\in\mathbb{R}$ and $|m|\geq5$, where $M\simeq1000$ are large enough.

One also easily obtain that
\begin{equation}\label{eq547}
  (|n_{i}|-6)^{-\frac{N_{i}}{2}}\lesssim\min_{x_{1}\in\widetilde{I}_{n_{i}}}\tilde{\chi}_{I_{0}}(x_{1}),
\end{equation}
for $i=1,2$ and every $|n_{1}|, |n_{2}|>15$, where $N_{1}, N_{2}\simeq1000$ are large enough.

Therefore, by using \eqref{eq545}, \eqref{eq546} and \eqref{eq547}, one can use the triangle inequality if $r\geq1$ and the subadditivity of $\|\cdot\|_{L^{r}}^{r}$ if $0<r<1$ to sum the contributions of $\overrightarrow{\Pi}_{a,\mathcal{R}}^{(n_{1},n_{2}),m}(f,g)$ together and obtain (we only present here the arguments for $r\geq1$, the cases $0<r<1$ can be treated similarly):
\begin{eqnarray}\label{eq548}
 &&\quad\,\|\sum_{|n_{1}|,|n_{2}|>15}\sum_{|m|\geq5}\overrightarrow{\Pi}_{a,\mathcal{R}}^{(n_{1},n_{2}),m}(f,g)\|_{L^{r}}\\
 \nonumber &\lesssim&\sum_{|n_{1}|,|n_{2}|>15}\sum_{|m|\geq5}[\prod_{i=1,2}\frac{1}{(|n_{i}|-6)^{N_{i}}}]\frac{1}{(|m|-2)^{M}}
 \|f\chi_{\widetilde{I}_{n_{1}}}\tilde{\chi}_{J_{m}}\|_{L^{p}}\|g\chi_{\widetilde{I}_{n_{2}}}\tilde{\chi}_{J_{m}}\|_{L^{q}}\\
 \nonumber &\lesssim&\sum_{|n_{1}|,|n_{2}|>15}\sum_{|m|\geq5}[\prod_{i=1,2}\frac{1}{(|n_{i}|-6)^{\frac{N_{i}}{2}}}]\frac{1}{(|m|-2)^{\frac{M}{2}}}
 \|f\tilde{\chi}_{R_{00}}\|_{L^{p}}\|g\tilde{\chi}_{R_{00}}\|_{L^{q}}\\
 \nonumber &\lesssim&\|f\tilde{\chi}_{R_{00}}\|_{L^{p}}\cdot\|g\tilde{\chi}_{R_{00}}\|_{L^{q}}.
\end{eqnarray}
Similar to the proof of estimate \eqref{eq548}, we can get the same estimates for the other different three cases, then we insert these estimates into the decomposition \eqref{eq542} and finally get the estimates of $\overrightarrow{\Pi}_{a,\mathcal{R},III}^{(2),(0,0,\vec{0})}$ as follows
\begin{eqnarray}\label{eq549}
 \|\overrightarrow{\Pi}_{a,\mathcal{R},III}^{(2),(0,0,\vec{0})}(f,g)\|_{L^{r}(\mathbb{R}^{2})}
   &=&\|\sum_{n_{1},n_{2}\in\mathbb{Z}}\sum_{|m|\geq5}\overrightarrow{\Pi}_{a,\mathcal{R}}^{(n_{1},n_{2}),m}(f,g)\|_{L^{r}(\mathbb{R}^{2})}\\
 \nonumber &\lesssim&\|f\tilde{\chi}_{R_{00}}\|_{L^{p}(\mathbb{R}^{2})}\cdot\|g\tilde{\chi}_{R_{00}}\|_{L^{q}(\mathbb{R}^{2})},
\end{eqnarray}
provided that $1<p, \, q\leq\infty$ and $\frac{1}{r}=\frac{1}{p}+\frac{1}{q}>0$, this concludes our estimates of the hybrid term $\overrightarrow{\Pi}_{a,\mathcal{R},III}^{(2),(0,0,\vec{0})}$.

As to the other hybrid term $\overrightarrow{\Pi}_{a,\mathcal{R},IV}^{(2),(0,0,\vec{0})}$, by symmetry, we can estimate it in a completely similar way as $\overrightarrow{\Pi}_{a,\mathcal{R},III}^{(2),(0,0,\vec{0})}$ by exchanging our arguments on variables $x_{1}$ and $x_{2}$, and finally obtain that
\begin{equation}\label{eq5410}
  \|\overrightarrow{\Pi}_{a,\mathcal{R},IV}^{(2),(0,0,\vec{0})}(f,g)\|_{L^{r}(\mathbb{R}^{2})}\lesssim
  \|f\tilde{\chi}_{R_{00}}\|_{L^{p}(\mathbb{R}^{2})}\cdot\|g\tilde{\chi}_{R_{00}}\|_{L^{q}(\mathbb{R}^{2})},
\end{equation}
provided that $1<p, \, q\leq\infty$ and $\frac{1}{r}=\frac{1}{p}+\frac{1}{q}>0$, this concludes our estimates of the hybrid terms $\overrightarrow{\Pi}_{a,\mathcal{R},III}^{(2),(0,0,\vec{0})}$ and $\overrightarrow{\Pi}_{a,\mathcal{R},IV}^{(2),(0,0,\vec{0})}$.

\subsection{Remarks on estimates for bilinear operators involved in decomposition \eqref{paraproducts decompostion} which contain at least one of components $\Pi^{i}_{ll}$ ($i=1,2$) in tensor products}

From the estimates of the standard discrete paraproduct operator corresponding to bilinear operator $\Pi^{1}_{lh}\otimes\Pi^{2}_{hl}$ presented in subsections 5.2-5.4, we realize when the supports of $f, \, g$ and dyadic rectangle $R$ are all close to $R_{00}$ in one direction (i.e. $I\subseteq5I_{0}$ or $J\subseteq5J_{0}$) but at least one of the supports of $f, \, g$ are far away from $R_{00}$ in the other direction, we need to apply the one-parameter paraproducts estimates (Theorem \ref{one-parameter paraproducts}) with respect to $x_{1}$ or $x_{2}$ variable, which is unfortunately inapplicable for the discretized operators $\overrightarrow{\Pi}_{a,\mathcal{R}}^{(2),(0,0,\vec{0})}$ corresponding to bilinear operators that contain at least one of $\Pi^{1}_{ll}$ or $\Pi^{2}_{ll}$ in the tensor products.

Indeed,
by a completely similar discretization procedure described in Section 4, one can reduce these seven bilinear operators
$$\Pi^{1}_{lh}\otimes\Pi^{2}_{ll}, \,\,\, \Pi^{1}_{hl}\otimes\Pi^{2}_{ll}, \,\,\, \Pi^{1}_{hh}\otimes\Pi^{2}_{ll}, \,\,\, \Pi^{1}_{ll}\otimes\Pi^{2}_{ll}, \,\,\, \Pi^{1}_{ll}\otimes\Pi^{2}_{hh}, \,\,\, \Pi^{1}_{ll}\otimes\Pi^{2}_{hl}, \,\,\, \Pi^{1}_{ll}\otimes\Pi^{2}_{lh}$$
appearing in the decomposition \eqref{paraproducts decompostion} of $T_{a}^{(2),\vec{0}}$ to averages of discrete bilinear paraproduct operators of the form \eqref{eq25} with restrictions $|I|,|J|\lesssim1$, and for at least one of the two dyadic interval families $\mathcal{I}$ and $\mathcal{J}$ (here we assume the tensor product contains $\Pi^{1}_{ll}$ and hence suppose it is dyadic interval family $\mathcal{I}$ without loss of generality), one has $|I|\sim1$ for every $I\in\mathcal{I}$ and at least two of the families of $L^{2}$-normalized bump functions $(\varphi_{I}^{j})_{I\in\mathcal{I}}$ for $j=1,2,3$ are nonlacunary. Therefore, different from the operator $T_{a,(lh,hl)}^{(2),\vec{0}}$, these seven operators can't be reduced to averages of classical discrete bilinear paraproduct operators of the form \eqref{eq25} which is applicable for Theorem \ref{paraproduct estimates}, even both the components $\Pi^{1}_{ll}$ and $\Pi^{2}_{ll}$ are inapplicable for Theorem \ref{one-parameter paraproducts}.

However, if the supports of both $f$ and $g$ are not far from the rectangle $R_{00}:=I_{0}\times J_{0}:=[-1,1]\times[-1,1]$, without loss of generality, we consider the operator $\Pi^{1}_{ll}\otimes\Pi^{2}_{ll}$, the cutoffs $f\chi_{15R_{00}}$ and $g\chi_{15R_{00}}$, since $\widetilde{\chi}_{R_{00}}$ is bounded from below on $15R_{00}$, we deduce from Coifman-Meyer theorem (Theorem \ref{22-multiplier}) that
\begin{equation}\label{eq46}
  \|\Pi^{1}_{ll}\otimes\Pi^{2}_{ll}(f\chi_{15R_{00}},g\chi_{15R_{00}})\cdot\varphi_{0}^{'}\otimes\varphi_{0}^{''}\|_{L^{r}(\mathbb{R}^{2})}\lesssim \|f\widetilde{\chi}_{R_{00}}\|_{L^{p}(\mathbb{R}^{2})}\cdot\|g\widetilde{\chi}_{R_{00}}\|_{L^{q}(\mathbb{R}^{2})}
\end{equation}
for any $1<p, q\leq\infty$ and $\frac{1}{r}=\frac{1}{p}+\frac{1}{q}>0$, which is acceptable for proving Proposition \ref{localized Coifman-Meyer}. Otherwise, if one of the supports of functions $f$, $g$ is far away from the rectangle $R_{00}$, note that $\Pi_{ll}^{i}$ are both summations of finite terms for $i=1,2$ and the dyadic intervals $I\in\mathcal{I}$ (corresponding to $\Pi^{1}_{ll}$), $J\in\mathcal{J}$ (corresponding to $\Pi^{2}_{ll}$) satisfy $|I|\sim1$ and $|J|\sim1$, so we don't need any decay factors or one-parameter paraproducts estimates (Theorem \ref{one-parameter paraproducts}, which can't be applied to $\Pi_{ll}^{1}$ and $\Pi_{ll}^{2}$) to make sure the summations $\sum_{I\subseteq5I_{0}}$ or $\sum_{J\subseteq5J_{0}}$ converge (see subsection 5.2 and 5.4), since both $\sum_{I\subseteq5I_{0}}$ and $\sum_{J\subseteq5J_{0}}$ are finite summations for $|I|\sim1$ and $|J|\sim1$. It is clear from the proof presented in subsections 5.1 - 5.4 that the other parts of our arguments have nothing to do with the properties whether the families of $L^{2}$-normalized bump functions $(\varphi_{I}^{j})_{I\in\mathcal{I}}$ and $(\varphi_{J}^{j})_{J\in\mathcal{J}}$ for $j=1,2,3$ are lacunary or not, we can deal with these seven operators in a quite similar way as $T_{a,(lh,hl)}^{(2),\vec{0}}$(see subsection 5.1, 5.2, 5.3 and 5.4).

\subsection{Conclusions}

By combining the estimate \eqref{eq529} for the main term $\overrightarrow{\Pi}_{a,\mathcal{R},I}^{(2),(0,0,\vec{0})}$, \eqref{eq536} for the error term $\overrightarrow{\Pi}_{a,\mathcal{R},II}^{(2),(0,0,\vec{0})}$, \eqref{eq549}, \eqref{eq5410} for the hybrid terms $\overrightarrow{\Pi}_{a,\mathcal{R},III}^{(2),(0,0,\vec{0})}$, $\overrightarrow{\Pi}_{a,\mathcal{R},IV}^{(2),(0,0,\vec{0})}$ and inserting them into the decomposition \eqref{eq511}, we finally obtain the estimates for the localized and discrete bilinear paraproduct operator $\overrightarrow{\Pi}_{a,\mathcal{R}}^{(2),(0,0,\vec{0})}$ as follows:
\begin{eqnarray}\label{eq551}
&&\quad\|\overrightarrow{\Pi}_{a,\mathcal{R}}^{(2),(0,0,\vec{0})}(f,g)\|_{L^{r}(\mathbb{R}^{2})}\\
\nonumber &\lesssim&\|\overrightarrow{\Pi}_{a,\mathcal{R},I}^{(2),(0,0,\vec{0})}(f,g)\|_{L^{r}(\mathbb{R}^{2})}
+\|\overrightarrow{\Pi}_{a,\mathcal{R},II}^{(2),(0,0,\vec{0})}(f,g)\|_{L^{r}(\mathbb{R}^{2})}
+\|\overrightarrow{\Pi}_{a,\mathcal{R},III}^{(2),(0,0,\vec{0})}(f,g)\|_{L^{r}(\mathbb{R}^{2})}\\
\nonumber&&\quad\quad\quad\quad\quad\quad\quad
+\|\overrightarrow{\Pi}_{a,\mathcal{R},IV}^{(2),(0,0,\vec{0})}(f,g)\|_{L^{r}(\mathbb{R}^{2})}
\lesssim\|f\tilde{\chi}_{R_{00}}\|_{L^{p}(\mathbb{R}^{2})}\cdot\|g\tilde{\chi}_{R_{00}}\|_{L^{q}(\mathbb{R}^{2})},
\end{eqnarray}
as long as $1<p, \, q\leq\infty$ and $\frac{1}{r}=\frac{1}{p}+\frac{1}{q}>0$, which completes the proof of Proposition \ref{paraproducts}.

This concludes the proof of our main result, Theorem \ref{main}.\\

{\bf Acknowledgements:}  The authors are greatly indebted to C. Muscalu who suggested to them this problem. They also wish to thank him very much
for his enlightening lectures on multi-linear and multi-parameter Coifman-Meyer theorems from which the authors have benefited, and for his comments on this paper.\\

\end{document}